\documentclass[11pt]{amsart}
\usepackage{amsfonts,amsmath,amssymb, latexsym}
\usepackage[dvips]{epsfig}
\usepackage{amscd}

\def\C{{\mathbb C}}
\def\cP{{\mathbb P}}
\def\R{{\mathbb R}}
\def\Q{{\mathbb Q}}
\def\Z{\mathbb Z}

\def\E{\mathbb E}

\def\CE{{\mathcal E}}
\def\CF{{\mathcal F}}

\def\CM{{\mathcal M}}
\def\CN{{\mathcal N}}
\def\CO{{\mathcal O}}
\def\CP{{\mathcal P}}

\newcommand{\oO}{\mathcal{O}}
\newcommand{\Coh}{\mathop{\rm Coh}\nolimits}
\newcommand{\rank}{\mathop{\rm rank}\nolimits}

\newcommand{\mM}{\mathcal{M}}
\newcommand{\ch}{\mathop{\rm ch}\nolimits}
\newcommand{\DT}{\mathop{\rm DT}\nolimits}
\newcommand{\eE}{\mathcal{E}}
\newcommand{\hH}{\mathcal{H}}
\newcommand{\dR}{\mathbf{R}}
\newcommand{\dL}{\mathbf{L}}
\newcommand{\zZ}{\mathcal{Z}}
\newcommand{\iI}{\mathcal{I}}
\newcommand{\NS}{\mathop{\rm NS}\nolimits}
\newcommand{\vV}{\mathcal{V}}
\newcommand{\aA}{\mathcal{A}}
\newcommand{\td}{\mathop{\rm td}\nolimits}

\def\b{\beta}

\newcommand{\Tr}{\mathrm{Tr}}
\newcommand{\Mbar}{\overline{\mathcal{M}} }

\newcommand{\tN}{\tilde{N}}
\newcommand{\tq}{\tilde{q}}

\newcommand{\XX}{X_0,X_1,X_2}

\newcommand{\z}{Z_0,Z_1}
\newcommand{\vir}{\mathrm{vir}}
\newcommand{\ev}{\mathrm{ev}}

\newcommand{\dbar}{\bar{\partial}}
\newcommand{\End}{\mathrm{End}}

\newcommand{\Hilb}{\mathrm{Hilb}}
\newcommand{\Sym}{\mathrm{Sym}}
\newcommand{\Mor}{\mathrm{Mor}}
\newcommand{\ad}{\mathrm{ad}}
\newcommand{\fg}{\mathfrak{g}}
\newcommand{\fsu}{\mathfrak{su}}
\newcommand{\fu}{\mathfrak{u}}

\newcommand{\tD}{\tilde{D}}

\newtheorem{fact}{Fact}
\newtheorem{theorem}{Theorem}

\newtheorem{conj}{Conjecture}
\newtheorem{lemma}{Lemma}

\newtheorem{remark}{Remark}

\begin{document}

\address{Sergei Gukov, Walter Burke Institute for Theoretical Physics, California Institute of Technology, Pasadena, CA 91125, USA}
\address{Sergei Gukov, Max-Planck-Institut f\:{u}r Mathematik, Vivatsgasse 7, D-53111 Bonn, Germany}
\email{gukov@theory.caltech.edu}
\address{Chiu-Chu Melissa Liu, Department of Mathematics, Columbia University, 2990 Broad- way, New York, NY 10027, USA}
\email{ccliu@math.columbia.edu}

\address{Artan Sheshmani, Center for Mathematical Sciences and Applications, Harvard University, 20 Garden Street, Cambridge, MA, 02138}
\email{artan@cmsa.fas.harvard.edu}
\address{Artan Sheshmani, Center for Quantum Geometry and Moduli spaces, Department of Mathematics, Aarhus University, Ny Munkegade 118, Building 1530, 8000 Aarhus C, Denmark}

\address{Shing-Tung Yau, Department of Mathematics, Harvard University, Cambridge, MA 02138, USA}
\email{yau@math.harvard.edu}
\title[Local theory of surfaces]{On topological approach to local theory of surfaces in Calabi-Yau threefolds}
\author[Gukov, Liu, Sheshmani, Yau]{Sergei Gukov, Chiu-Chu Melissa Liu, Artan Sheshmani and Shing-Tung Yau}
\date{\today}
\maketitle

\begin{abstract}
We study the web of dualities relating various enumerative invariants,
notably Gromov-Witten invariants and invariants that arise in topological
gauge theory.
In particular, we study Donaldson-Thomas gauge theory and its reductions
to $D=4$ and $D=2$ which are relevant to the local theory of surfaces
in Calabi-Yau threefolds.
\end{abstract}

{\small \tableofcontents }


\section{Introduction}

Not so long ago, enumerative invariants were completely unfamiliar to physicists and
were regarded more as a curiosity rather than serious mathematics in the math world.
At present, enumerative invariants are heavily used in theoretical physics
and play a major role in many branches of pure mathematics:
symplectic geometry, algebraic geometry, mirror symmetry, and gauge theory.

In gauge theory, the prominent examples of enumerative invariants are Donaldson polynomials
and Seiberg-Witten invariants, which help to distinguish different smooth structures on 4-manifolds.
In recent years, other 4-manifold invariants have been introduced by changing the gauge theory
(the PDE's and the ``counting'' problem) or, by changing dimension, similar gauge theory invariants
were defined on higher-dimensional manifolds.
Notable examples include Donaldson-Thomas (DT) invariants which, roughly speaking, can be viewed as six-dimensional
analogues of Donaldson invariants on 4-manifolds~\cite{DT}.

In symplectic geometry and mirror symmetry, enumerative invariants play an equally important role;
it is almost impossible to imagine these fields without enumerative invariants.
Most of the original versions of such invariants were based on the problem of curve ``counting'' of some form.
More recently, however, this was generalized to a much larger framework where the role of curves is less special
and, in the same time, the borderline between gauge theory and symplectic/algebraic geometry started to blur.
For example, Donaldson-Thomas invariants that we mentioned earlier and that originally were introduced via gauge theory,
nowadays are actually more widely regarded as part of algebraic geometry and mirror symmetry, since they are defined as the intersection numbers over the moduli space of sheaves (say for example the ideal sheaves of one dimensional subschemes of the ambient variety).

String theory also stimulated the process of erasing borders.
It led to many dualities which relate enumerative invariants in different dimensions and even across the fields.
For example, gauge theory itself appears in string theory through {\it open} strings,
whereas curve counting (Gromov-Witten theory) is naturally a part of the {\it closed} string theory.
Therefore, string dualities which relate open and closed strings (also known as gauge/gravity dualities)
provide a fundamental, conceptual reason why one might hope to have a relation
between gauge theory enumerative invariants and, say, Gromov-Witten invariants.

In this paper, we explore a rich web of inter-relations
and dualities between different enumerative invariants.

\subsection*{Acknowledgements} We would like to thank
Jun Li, Davesh Maulik, and Richard Thomas for helpful conversations.
The work of S. Gukov is supported in part by the U.S. Department of Energy, Office of Science,
Office of High Energy Physics, under Award Number DE-SC0011632 and in part
by the ERC Starting Grant no. 335739 ``Quantum fields and knot homologies''
funded by the European Research Council under the European Union Seventh Framework Programme.
The work of C.-C. Liu is partially supported by NSF grants DMS-1206667 and DMS-1159416.
A. Sheshmani would like to thank Kavli IPMU, MIT, Harvard and the Institute Henri Poincar\'{e} (IHP)  for creating the opportunity of initiating the discussions about the current article. The work of A. Sheshmani was partially supported by the World Premier International Research Center Initiative (WPI Initiative), MEXT, Japan. The work of S.-T. Yau is partially supported by NSF grants DMS 1308244, DMS-159412, PHY-1306313, and PHY-0937443.


\section{Preliminaries}

\subsection{Gromov-Witten theory} \label{sec:reviewGW}

\subsubsection{Gromov-Witten invariants
and Taubes's Gromov invariants of symplectic 4-manifolds}\label{sec:severi}

In this subsection, we briefly review certain Gromov-Witten invariants
of compact symplectic 4-manifolds and Taubes's Gromov invariants
\cite{Taubesbook,IP}.

Let $(S,\omega)$ be a compact symplectic 4-manifold with
an almost complex structure $J$ compatible with the symplectic
structure $\omega$, so that $(S,\omega,J)$ is an almost K\"{a}hler manifold.
If $J$ is integrable then $(S,\omega,J)$ is a K\"{a}hler manifold.

Let $\Mbar_{g,n}(S,\beta)$ be the moduli space
of genus $g$, degree $\beta\in H_2(S;\Z)$ stable maps
to $(S,\omega,J)$. The virtual dimension of $\Mbar_{g,n}(S,\beta)$ is
given by
$$
\int_{\beta} c_1(S) +g-1+n = -K_S\cdot \beta + g-1 + n,
$$
where $K_S$ is the canonical class of $S$. 
 
Let $\nu\in H^4(S;\Z)$ be the Poicare dual of
the class of a point. We define
\begin{equation}\label{eqn:NS}
N^S_{g,\beta}:=\left\{ \begin{array}{ll}
\int_{[\Mbar_{g,n}(S,\beta)]^{\vir} } \prod_{i=1}^n \ev_i^*(\nu), &
\textup{ if } n= -K_S\cdot \beta+g-1\geq 0,\\
0, & \textup{ if }-K_S\cdot \beta+g-1 <0,
\end{array}\right.
\end{equation}
where $\ev_i:\Mbar_{g,n}(S,\beta)\to S$ is the evaluation at
the $i$-th marked point. By the divisor equation, any genus $g$, degree
$\beta$ primary Gromov-Witten invariants of $S$ can be reduced to
$N^S_{g,\beta}$. In general, Gromov-Witten invariants count
connected parametrized curves, and are rational numbers instead of integers. 

By adjunction formula, the genus of a connected, embedded smooth
$J$-holomorphic curve
in class $\beta\in H_2(S;\Z)$ is given by
$$
1+\frac{1}{2}(\beta\cdot \beta + K_S\cdot \beta).
$$
Gromov invariant $N_{g,\beta}^S$ counts connected embedded curves
of genus $g$ in class $\beta$ passing through $n$
generic points if the following three conditions holds:
\begin{enumerate}
\item $g=1+ \frac{1}{2}(\beta\cdot \beta +K_S \cdot \beta)$.
\item $n=\frac{1}{2}(\beta\cdot\beta -K_S \cdot \beta)\geq 0$.
\item $\beta \notin T$, where $T=\{ \beta\in H_2(S,\Z)\mid \beta\cdot K_S=\beta\cdot \beta=0\}$
is the set of {\em toroidal casses}.
\end{enumerate}
When the above conditions hold, the moduli space of $J$-holomorphic curves of genus $g$ in class $\beta$ passing through
$n$ points contains no multiply covered curves for generic $J$ (see Ruan \cite{R}, Taubes \cite{Taubesbook}).
By \cite[Section 4]{IP},
$$
N_{g,\beta}^S = \widehat{Gr}_S(\beta)
$$
where $\widehat{Gr}_S(\beta)$ is the connected version of Taubes's Gromov invariant in class $\beta$.
Taube's Gromov invariant $Gr_S(\beta)$, which counts possibly disconnected curves, is equivalent to 
certain Seiberg-Witten invariant by Taubes's famous theorem. 

In this paper, we will assume $J$ is integrable, so that
the pair $(\omega,J)$ is a K\"{a}hler structure on $S$.

\subsubsection{Non-equivariant local Gromov-Witten invariants of surfaces in a Calabi-Yau threefold}\label{sec:local}
Let $S$ be a projective Fano surface, and let $X_S$ be the total
space of the canonical line bundle $K_S$ over $S$. Then $X_S$ is a noncompact Calabi-Yau 3-fold. Given
any nonzero $\beta\in H_2(S;\Z)\cong H_2(X_S;\Z)$, we have
$$
\Mbar_{g,0}(X_S,\beta)=\Mbar_{g,0}(S,\beta)
$$
as Deligne-Mumford stacks. The virtual dimensions of the perfect
obstruction theories on $\Mbar_{g,0}(X_S,\beta)$ and on
$\Mbar_{g,0}(S,\beta)$ are $0$ and $-K_S\cdot \beta +g-1$,
respectively, and
$$
N^{X_S}_{g,\beta}=\int_{[\Mbar_{g,0}(X_S,\beta)]^\vir} 1
=\int_{[\Mbar_{g,0}(S,\beta)]^\vir} e(V_{g,\beta})
$$
is the genus $g$, degree $\beta$ Gromov-Witten invariant of $X_S$,
where $V_{g,\beta}$ is a rank $-K_S\cdot\beta+g-1$ vector bundle
over $\Mbar_{g,0}(S,\beta)$ whose fiber over $[f:C\to S]$ is
$H^1(C,f^* K_S)$. Note that because $\beta$ is nonzero and $-K_S$
is ample, we have $\deg(f^* K_S)<0$, so $H^0(C, f^* K_S)=0$.
(See \cite[Section 5.2]{CKYZ}.) The invariants $N^{X_S}_{g,\beta}$ can be viewed
as local invariants of the Fano surfaces $S$ in a Calabi-Yau 3-fold.

When $S$ is projective but not Fano, the Gromov-Witten invariants of $X_S$ are
not necessarily defined, due to the non-compactness of
moduli spaces $\Mbar_{g,0}(X_S,\beta)$. Observe that, when
the class $\beta$ satisfies $-K_S\cdot \beta >0$,
$V_{g,\beta}$ is a rank $-K_S\cdot \beta + g-1$ vector bundle
over $\Mbar_{g,0}(S,\beta)$. This allows us to define the genus $g$, degree $\beta$
Gromov-Witten invariant of $X_S$ by
\begin{equation}\label{eqn:non-equivariant-local}
N_{g,\beta}^{X_S}:= \int_{[\Mbar_{g,0}(S,\beta)]^\vir} e(V_{g,\beta}).
\end{equation}

Note that $e(V_{g,\beta})=1$
if and only if $H^1(C,f^* K_S)=0$ for all $[f:C\to X_S]\in
\Mbar_{g,0}(S,\beta)$. In this case
$$
\int_{[\Mbar_{g,0}(X_S,\beta)]^\vir}1=
\int_{[\Mbar_{g,0}(S,\beta)]^\vir}1
$$
is also the genus $g$, degree $\beta$ Gromov-Witten invariant of
the surface $S$. We have
$$
-K_S\cdot \beta + g-1 =0
$$
where $-K_S\cdot \beta >0$, so we must have $g=0$ and $K_S\cdot
\beta =-1$.  We summarize the above discussion in the following lemma.
\begin{lemma}
Let $S$ be a projective surface, and suppose $\beta\in H_2(S,\Z)$ satisfies
$K_S\cdot \beta=-1$. Then
$$
N^{X_S}_{0,\beta}=N^S_{0,\beta}
$$ 
where the left hand side is defined by \eqref{eqn:non-equivariant-local} and
the right hand side is defined as in Section \ref{sec:severi}.   
\end{lemma}

In the remainder of this subsubsection, we assume 
$$
S \textup{ is Fano},\quad g=0,\quad K_S\cdot\beta = -1.
$$
Then a curve in class $\beta$ cannot be multiply covered.
We claim that there is no bubbling. Indeed if
$\beta=\beta_1+\beta_2$ where $\beta_1,\beta_2$ are nonzero
effective class, we have $-K_S\cdot \beta_i>0$ for $i=1,2$
so $-K_S\cdot\beta \geq 2$, which is a contradiction.
Similarly, $\beta$ cannot be represented by
a disconnected holomorphic curve. Combining the above discussion 
and the discussion Section \ref{sec:severi}, we conclude: 
\begin{lemma}
Let $S$ be a Fano surface, and suppose that $\beta\in H_2(S,\Z)$
satsifies
$$
K_S\cdot \beta = \beta\cdot \beta =-1.
$$
Then
$$
N_{0,\beta}^{X_S} = N_{0,\beta}^S = \widehat{Gr}_S(\beta) = Gr_S(\beta).
$$
\end{lemma}

Let $B_k$ denote $\cP^2$ blowup at $k$ generic points. By classification of surfaces,
Fano surfaces are: $\cP^2$, $\cP^1\times \cP^1$, and del Pezzo surfaces $B_k$ where $1\leq k\leq 8$.
If $S=\cP^2$ or $\cP^1\times \cP^1$ then there is no class $\beta$ satisfying $K_S\cdot \beta =-1$.
If $S=B_k$ and $e_1,\ldots, e_k \in H_2(S,\Z)$ are classes of exceptional divisors, then 
$K_S\cdot e_i = e_i\cdot e_i=-1$.

\subsubsection{Equivariant local Gromov-Witten invariants of surfaces in a Calabi-Yau threefold}

It is possible to define {\em residue invariants} in the presence of
a torus action \cite[Section 2.1]{BP}. We first consider the special case of 
toric surfaces. We have
$$
N^{X_S}_{g,\beta}=\int_{[\Mbar_{g,0}(X_S,\beta)]^\vir} 1
=\int_{[\Mbar_{g,0}(X_S,\beta)^T]^\vir }
\frac{1}{e_T(N^\vir)}.
$$
where $T$ acts on $X_S$, and $e_T(N^\vir)$ is the $T$-equivariant
Euler class of the virtual normal bundle of $\Mbar_{g,0}(X_S,\beta)^T$
in $\Mbar_{g,0}(X_S,\beta)$. When the $T$-action on $X_S$ induced from
a $T$-action on the base $S$, so that the $T$-action on $X_S$ acts trivially
on the canonical line bundle of the toric threefold $X_S$, these invariants
can be evaluated by the algorithm of the topological vertex 
\cite{AKMV}. Indeed, for local toric surfaces, only 2-leg vertex is needed;
the algorithm in this case is described in \cite{AMV}. The topological
vertex is proved in the 1-leg case in \cite{LLZ1, OP} and the 2-leg case
in \cite{LLZ2, LLLZ}. The full 3-leg case is a consequence of
the proof of the GW/DT corresspondence of toric 3-folds in \cite{MOOP}.

In general, if there is a torus action
on $X_S$ such that $X_S^T$ is contained in the zero section $S$
and $\Mbar_{g,0}(X_S,\beta)^T = \Mbar_{g,0}(S,\beta)^T$ is compact, 
we define
\begin{equation}\label{eqn:NXt}
N^{X_S}_{g,\beta}(t_1,\ldots,t_r)=\int_{[\Mbar_{g,0}(X_S,\beta)^T]^\vir }
\frac{1}{e_T(N^\vir)} \in \Q(t_1,\ldots,t_r)
\end{equation}
where $t_1,\ldots,t_r$ are the generators of
the $T=(\C^*)^r$ equivariant cohomology of a point. Alternatively,
\begin{equation}\label{eqn:NXtS}
N^{X_S}_{g,\beta}(t_1,\ldots,t_r)=\int_{[\Mbar_{g,0}(S,\beta)^T]^\vir }
\frac{e_T(H^1)/e_T(H^0)}{e_T(N^\vir)} \in \Q(t_1,\ldots,t_r)
\end{equation}
where $e_T(N^\vir)$ is the $T$-equivariant Euler class of the virtual
normal bundle of $\Mbar_{g,0}(S,\beta)^T$ in $\Mbar_{g,0}(S,\beta)$.
Here $H^1$ and $H^0$ are vector bundles on each connected component
of $\Mbar_{g,0}(S,\beta)^T$, and their fibres over $[f:C\to S]$
are $H^1(C,f^* K_S)$ and $H^0(C,f^* K_S)$, respectively.


\subsection{Donaldson-Thomas theory} \label{sec:reviewDT}

We consider a general gauge group $G$ which can be any
compact connected Lie group. The complexification $G^\C$ of $G$ is
a connected reductive complex algebraic group. For example, when
$G=U(N)$, we have $G^\C=GL(N,\C)$.
The space of connections on a principal $G$-bundle $P$ over
a {\em K\"{a}hler} manifold can be identified with
the space of $\dbar$-operators on the principal $G^\C$ bundle
$P^\C =P\times_G G^\C$.


The Donaldson-Thomas theory is a topological gauge theory
defined on a complex Kahler 3-fold $X$ which, roughly speaking,
``counts'' solutions to the Hermitian Yang-Mills (HYM) equations
(see e.g. \cite[Definition 3.2]{AB}):
\begin{equation}\label{eqn:HYM}
F^{0,2}_A =0,\quad \Tr_{k_0} F_A^{1,1}=\tau.
\end{equation}
where
$A \in \Omega^1_{X}(\ad P)$ is the gauge connection,
$k_0$ is the K\"{a}hler form on $X$,
and $\tau$ is an element in the center of the Lie algebra $\fg$ of $G$.
For example, $\tau=0$ if $G$ is semisimple, and $\tau$ is a multiple of the
identity element if $G=U(N)$.

This six-dimensional TQFT can be obtained by a topological twist
of the dimensional reduction of the ten-dimensional
super-Yang-Mills \cite{BKSinger,ALSpence,BT,DT,BSV,INOV}.
Under the reduction,
$$
SO(10) \to SO(6) \times SO(4)
$$
and bosonic/fermionic fields decompose as
\begin{eqnarray}
{\rm bosons:} && {\bf 10} \to ({\bf 6},{\bf 1}) \oplus ({\bf 1},{\bf 4}) \cr
{\rm fermions:} && {\bf 16} \to ({\bf 4},{\bf 2}) \oplus (\bar {\bf 4}, \bar {\bf 2})
\end{eqnarray}
In addition, if $X$ is K\"{a}hler, we have $SO(6) \to U(3) \cong SU(3) \times U(1)$.
The twisted six-dimensional gauge theory is obtained by mixing the last $U(1)$
factor with $U(1) \subset SO(4)$. As a result, one finds a theory with $\CN_T = 2$
topological supersymmetries (nilpotent BRST symmetries) and the following field content:

\medskip
\centerline{\vbox{\halign{\quad # & \quad # & \quad # \cr
bosonic: & 1-form $A$ & \cr
& $(3,0)$-form $\varphi$ & \cr
& complex scalar $\phi$, $\bar \phi$ & \cr
fermionic: & 1-form $\psi^{1,0}$, $\psi^{0,1}$ & 3+3 \cr
& 2-form $\chi^{2,0}$, $\chi^{0,2}$ & 3+3 \cr
& 3-form $\psi^{3,0}$, $\psi^{0,3}$ & 1+1 \cr
& two scalars $\eta$, $\bar \eta$ & 1+1 \cr }}}

\medskip
\noindent
where all fields transform in the adjoint representation of the gauge
group $G$ and in the last column we counted the number of real components.
Notice that, altogether, the fermionic fields contain 16 real components.
If $X$ is a Calabi-Yau space, the resulting theory has $\CN_T = 4$
topological supersymmetries.

The partition function of the six-dimensional topological gauge theory
localizes on the solutions to the following equations
\begin{equation}\label{eqsprime}
\begin{aligned}
& F_A^{0,2}= \dbar_A^{\dagger} \bar{\varphi} \\
& \Tr_{k_0} F_A^{1,1} + *[\varphi,\bar{\varphi}] =\tau \\
& d_A \phi =0
\end{aligned}
\end{equation}
where the first two equations can be recognized as perturbations
of the HYM equations \eqref{eqn:HYM}.
The Bianchi identity and the first equation imply
$\dbar_A \dbar_A^{\dagger}\bar{\varphi}=0$, so when $X$ is compact,
the first equation is decoupled into two equations
$\dbar_A^{\dagger}\bar{\varphi}=0$ and $F_A^{0,2}=0$, and the
path integral localizes on the solutions to
\begin{equation}\label{eqn:extended}
\dbar_A^{\dagger} \bar{\varphi}=0, \quad F_A^{0,2}=0,\quad
\Tr_{k_0} F_A^{1,1} + *[\varphi,\bar{\varphi}] =\tau
\end{equation}
Note that the solutions to the HYM equations
\eqref{eqn:HYM} are also solutions to \eqref{eqn:extended}.
In \cite{HP}, the moduli space defined by \eqref{eqn:extended} is called the {\em extended}
moduli space of Einstein-Hermtian (i.e. HYM) connections. We have
$$
\CM_{\mathrm{HYM}}(X,G)\subset \CM_{\mathrm{DT}}(X,G),
$$
where $\CM_{\mathrm{HYM}}(X,G)$ and $\CM_{\mathrm{DT}}(X,G)$
are moduli spaces of solutions to \eqref{eqn:HYM} and \eqref{eqn:extended},
repectively. Given a holomorphic structure, the space of solutions
to $\dbar^{\dagger}_A \bar{\varphi}=0$ can be identified with
$H^0(X,K_X \otimes \ad(P))= H^3(X,\ad(P))^\vee$.

When $X$ is Calabi-Yau, $\varphi$ can be viewed as an element
of $H^0(X,\ad(P))$. If $P$ is stable, then $\varphi$ lies
in the center $Z(\fg)$ of $\fg$, so $[\varphi,\bar{\varphi}]=0$,
so the second and third equations in \eqref{eqn:extended} become
the HYM equations \eqref{eqn:HYM}.
 For example, $Z(\fsu(N))$ is trivial, and $Z(\fu(N))$ is spanned by
the identity matrix. When $X$ is Fano, i.e., $-K_X>0$, we have
$H^0(X,K_X\otimes \ad(P))=H^3(X,\ad(P))^\vee=0$.

Let $\CM$ denote the moduli space of gauge equivalence classes of solutions
to \eqref{eqsprime}. The tangent space of $\CM$ at a point $(A,\varphi)$
can be identified with $H^1$ of the following deformation
complex\footnote{Using $\dbar_A^2 \eta = [F_A^{0,2},\eta]$
and $F_A^{0,2}= \dbar_A^{\dagger} \bar{\varphi}$, one can check that
this is actually a complex:
$$
D_2\circ D_1(\eta)= D_2 (\dbar_A\eta, [\bar{\varphi},\eta])
=\dbar_A^2 \eta -\dbar_A^\dagger [\bar{\varphi},\eta]
=[\dbar_A^\dagger \bar{\varphi}, \eta]-\dbar_A^\dagger[\bar{\varphi},\eta]=0
$$}
\begin{equation} \label{sixdcplx}
0\to \Omega_{X}^{0,0}(\ad P_\C)\stackrel{D_1}{\to}
\Omega_{X}^{0,1}(\ad P_\C)\oplus \Omega_{X}^{0,3}(\ad P_\C)
\stackrel{D_2}{\to}\Omega^{0,2}_{X} (\ad P_\C )\to 0.
\end{equation}
When $X$ is Calabi-Yau, the index of the above deformation
complex is zero. We have
\begin{equation}
\dim T_{(A,\varphi)}\CM = \dim \CN_{(A,\varphi)}
\end{equation}
where $\CN_{(A,\varphi)}$ is $H^2$ of the above deformation
complex. In this case, from 1-loop determinants one finds that the partition
function computes the Euler characteristic of the obstruction bundle
$\CN$. We have
\begin{eqnarray*}
Z_{DT}^G &=& \int  DA D\varphi D(\ldots)  \exp(S(A,\varphi,\ldots))\\
&=&\int_{\CM}e(\CN)  \exp \int_X \Tr \big(
{g_s \over 3!} F \wedge F \wedge F
+ {1 \over 2} F \wedge F \wedge k_0
+ F \wedge \vartheta \big)
\end{eqnarray*}
where  $e(\CN)$ is the Euler class of $\CN$.
The integral over $X$ is topological in the sense that it is
constant on each connected component of $\CM$.
For example, when $G=U(N)$, let $\CM_{ch_1,ch_2,ch_3}$ be the
union of connected components of $\CM$ with fixed
Chern characters $ch_1,ch_2,ch_3$. Then
\begin{equation}\label{ZDT}
Z_{DT}^{U(N)}(q,t,T)=\sum_{ch_1,ch_2,ch_3}
\int_{\CM_{ch_1,ch_2,ch_3} } e(\CN)  q^{ch_3 (\CE)} t^{ch_2 (\CE)} T^{ch_1 (\CE)}
\end{equation}
where $q = e^{g_s}$, $t$, and $T$ are formal variables
that encode dependence of the Donaldson-Thomas partition
function on the Chern classes. More generally, when $X$ is not a Calabi-Yau manifold,
the index of the deformation complex \eqref{sixdcplx} is
\begin{equation}\label{eqn:index}
\int_X c_1(T_X)\left(ch_0ch_2-\frac{ch_1^2}{2} +\frac{ch_0^2-1}{24}c_2(T_X)\right)
=\mathrm{rank}_\C \CN -\dim_\C \CM
\end{equation}
and the correlation functions are of the form
\begin{equation}\label{eqn:correlation}
\int_{\CM_{ch_1,ch_2,ch_3}}e(\CN) \CO_1 \cdots \CO_k
\end{equation}
where $\CO_i\in H^*(\CM_{ch_1,ch_2,ch_3})$ are topological observables,
invariant under the BRST symmetry.

In order to define \eqref{eqn:correlation} mathematically,
one needs to compactify the moduli space $\CM$. By Uhlenbeck
and Yau's result \cite{UY}, the moduli space
of gauge equivalence classes of HYM $U(N)$-connections
can be identified with the moduli space of isomorphism
classes of polystable vector bundles of rank $N$. The later can be compatified
by the moduli space of rank $N$ semistable sheaves. In the algebro-geometric
framework, the deformation complex \eqref{eqsprime} is replaced
by the tangent-obstruction complex which encodes the deformation
theory of the moduli problem. When higher obstruction groups vanish,
one obtains virtual fundamental cycle $[\CM]^\vir$, and the mathematical
definition of the integral \eqref{eqn:correlation} is given by
\begin{equation}
\int_{[\CM_{ch_1,ch_2,ch_3}]^\vir} \CO_1 \cdots \CO_k\in \Z.
\end{equation}
The deformation theory and the virtual fundamental cycle
have been established for smooth projective Calabi-Yau and Fano 3-folds
in \cite{Th}.
In string theory, the Donaldson-Thomas theory with gauge
group $G=U(N)$ can be interpreted as a gauge theory on
$N$ D6-branes wrapped on $X$. In this interpretation,
the integer coefficients of the partition function $Z_{DT}$
count BPS states with D0, D2, D4, and D6-brane changes.
We remind that D-brane charges are described by the Mukai vector,
$Q \in H^0 (X) \oplus H^2 (X) \oplus H^4 (X) \oplus H^6 (X)$, given by
$$
Q = ch (\CE) \sqrt{Td(X)}
$$
These invariants are in fact the generalized DT invariants which were constructed and computed by the wall crossing techniques of Kontsevich-Soibelman \cite{KS} and Joyce-Song \cite{JS}. In connection with Gromov-Witten invariants which count holomorphic curves,
one considers moduli space of ideal sheaves of curves.
These are torsion free, rank 1 sheaves with trivial determinant, so they
correspond to abelian theory with $ch_1=0$.
The conjectural GW/DT correspondence \cite{MNOP1, MNOP2}
(proved for toric 3-folds \cite{MOOP} and local curves \cite{BP, OP2}) 
states that
$$
Z_{DT}'(q,t) = Z_{GW}'(g_s,t)
$$
where $Z_{DT}'(q,t)$ and $Z_{GW}'(g_s,t)$ are
the {\em reduced} DT and GW partition functions. Roughly speaking,
$Z_{DT}'(q,t)$ is obtained from $Z_{DT}(q,t)$, by removing the contribution
of point-like instantons, and $Z_{GW}'(g_s,t)$ is obtained from $Z_{GW}(g_s,t)$
by removing the contribution of constant maps.

In certain cases, the Donaldson-Thomas theory on $X$ can be
related to topological theories in lower dimensions.
In this paper, we describe two such relations which, in a sense,
are analogues of the corresponding relations
in four-dimensional topological gauge theory
studied, respectively, in \cite{BJSV} and \cite{MNVW}.

\subsection{Seiberg-Witten theory}
In this subsection, we briefly review the definition
of the Seiberg-Witten invariants \cite{Wmonopoles}
\begin{equation}
SW_S ~:~~ Spin^c (S) \to \Z
\end{equation}
and their relation \cite{Taubesbook} with the enumerative
invariants of $S$.

Up to the finite group $H^1 (S,\Z_2)$,
the set of Spin$^c$ structures on a 4-manifold $S$ is
parameterized by integral cohomology classes\footnote{Sometimes,
we parametrize the set of Spin$^c$ structures by $\lambda = {c
\over 2} \in H^2 (S,\Z) + w_2 (S)/2$, which is widely used in the
physics literature.} which reduce to $w_2 (X)$ mod 2
$$
{\rm Spin}^c (S) = \{ c \in H^2 (S,\Z) ~\vert~ c \equiv w_2 (S)
~{\rm mod}~2 \}
$$
Given a Spin$^c$ structure $c \in {\rm Spin}^c (S)$, let $L$ be
the corresponding Hermitian line bundle, and $S^{\pm}_L$ the
corresponding spinor bundles. The Seiberg-Witten monopole
equations are equations for a pair $(A,\Psi)$, where $A$ is a
unitary connection on $L$ and $\Psi$ is a smooth section of
$S^+_L$. In order to write the equations, we need to introduice
the Dirac operator
$$
D_A ~:~~ S^+_L \to S^-_L
$$
and a map
$$
\Omega^0 (S^+_L) \to \Omega^0 ( {\rm ad}_0 ~ S^+_L )
$$
$$
\Psi \mapsto i (\Psi \otimes \Psi^*)_0
$$
where $\Psi^*$ is the adjoint of $\Psi$ and ${\rm ad}_0 ~ S^+_L$
is the subbundle of the adjoint bundle of $S^+_L$ consisting of
the traceless skew-Hermitian endomorphisms. Then, the
Seiberg-Witten equations take the form \cite{Wmonopoles}
$$
D_A \Psi = 0
$$
$$
F^+_A - i (\Psi \otimes \Psi^*)_0 = 0
$$
where in the second equation we used the identification between
the space of self-dual 2-forms and skew-Hermitian automorphisms of
the positive spin representation, ${\rm ad}_0 ~ S^+_L \cong
\Lambda ^2_+ $.

Let $\CM_c$ be the moduli space of solutions to the Seiberg-Witten
equations. It has virtual dimension
$$
d_c = {1 \over 4} (c^2 - 2 \chi (S) - 3 \sigma (S) )
$$
The Seiberg-Witten invariants are defined as\footnote{Even though
some of our examples below include 4-manifolds with $b_1 > 0$,
we will not consider classes in $H^1 (S;\R)$. In particular,
we will not worry about the orientation of $H^1 (S;\R)$.}
$$
SW_S (c) = \int_{\CM_c} a_D^{d_c/2}
$$
where $a_D$ is a 2-form which represents the first Chern class of
the universal line bundle on the moduli space $\CM_c$. Notice,
that the dimension of the Seiberg-Witten moduli space and $b_2^+ +
b_1$ have opposite parity. This fact can be used to establish
vanishing of Seiberg-Witten invariants when $b_2^+ + b_1$ is even.
Another useful fact is a vanishing theorem \cite{Wmonopoles},
which says that $SW_S (c)=0$ whenever $S$ admits a metric of
positive scalar curvature. This theorem is especially helpful
in some of our examples, where the existence of a positive
curvature K\"{a}hler-Einstein metric on $S$ is directly
related to the existence of a complete Ricci flat K\"{a}hler
metric (Calabi-Yau metric) on $X_S$ \cite{Y,TY1, TY2}.

Among the examples that we consider, there are 4-manifolds with
$b_2^+ =1$. In such cases, Seiberg-Witten invariants are not quite
topological invariants; they are only piece-wise constant and
exhibit discontinuous behavior under wall crossing that we briefly
review below.

Let $H$ be the hyperbolic space
$$
H = \{ h \in H^2_{DR} (S) ~\vert~ h^2=1 \}
$$
This space has two connected components and the choice of one of
them orients the lines $H^{2,+}_g (S)$ for all metrics $g$. Having
fixed a component $H_0$ of $H$, every metric defines a unique
self-dual 2-form $\omega$ of length 1 with $[\omega] \in H_0$,
called the period point, $\omega^2 = 1$.
We have $H^{2,+} (S;\R) \cong \R \omega$ and $c$ can be written
as $c = c_+ + c_-$ where $c_+ = (c \cdot \omega) \omega$.


The space $H$ is divided into a set of chambers $C$, and for every
metric $g$ whose self-dual 2-form $\omega$ lies in $\R_{+} \cdot C$
one can define the corresponding Seiberg-Witten invariants $SW_{S,g} = SW_{S,C}$.
Notice, that changing $\omega$ to $- \omega$ corresponds to changing
the orientation of the Seiberg-Witten moduli space, so that
\begin{equation}
SW_{S,-C} = - SW_{S,C}
\label{swreverse}
\end{equation}
The Seiberg-Witten invariants $SW_{S,C}$ can be defined
using the perturbed (or ``twisted'') monopole
equations\footnote{In topological gauge theory, the perturbation
by the self-dual 2-form $\beta^+$ corresponds to perturbation by
a mass term \cite{WKahler}.}
$$
D_A \Psi = 0
$$
$$
\left( F_A  + 2 \pi i \beta \right)^+ = i (\Psi \otimes \Psi^*)_0
$$
where $\beta$ is a closed 2-form with $b = [\beta] \in H^2 (S;\R)$.
Equivalently, we consider a Spin$^c$ structure given by $(c-b)$.
Each characteristic class $c \in H^2 (S,\Z)$, that is an integral
lift of $w_2 (S)$, defines four chambers of type $c$:
$$
C_{H_0 , \pm} = \{ (h,b) \in H_0 \times H^2_{DR} (S) ~\vert~ \pm (c-b) \cdot h < 0 \}
$$
where $H_0$ is one of the components of $H$. Following the
standard conventions, we denote the corresponding Seiberg-Witten
invariants by $SW^{\pm}_{S,H_0} (c)$ (or simply by $SW^{\pm}_{S} (c)$).
Notice, according to (\ref{swreverse}), we have $SW_{S,-H_0}^{\pm} (c) = - SW_{S,H_0}^{\mp} (c)$.

The difference between the Seiberg-Witten invariants across a wall
is given by a simple formula, which for a simply-connected surface $S$ takes the form
$$
SW_S^+ (c) - SW_S^- (c) = 1
$$
For example, for $S=\cP^2$, we have Spin$^c( \cP^2 ) \cong (2\Z +
1)h$, where $h$ denotes the first Chern class of $\CO(1)$,
$h^2=1$. And for every $c \equiv h$ mod 2, there are two chambers
of type $c$:
$$
C_{H_0,\pm} = \{ (h,b) \in H_0 \times H_{DR}^2 (S) ~\vert~ \pm
(c-b) \cdot h < 0 \}
$$
where $H_0= \{ h \}$ and $d_c = {1 \over 4} (c^2-9)$. The
corresponding Seiberg-Witten invariants are
\begin{eqnarray}
&& SW_{S}^+ (c) = 1 \quad\quad c \cdot h \ge 3 \cr && SW_{S}^+ (c)
= 0 \quad\quad c \cdot h <3 \cr && SW_{S}^- (c) = -1 \quad\quad c
\cdot h \le -3 \cr && SW_{S}^- (c) = 0 \quad\quad c \cdot h
> -3
\end{eqnarray}
More generally, for K\"ahler surfaces with $b_2^+ = 1$ and $b_1 = 0$,
we have $SW_{S}^{\pm} (c) = \{ 0,1 \}$ or $SW_{S}^{\pm} (c) = \{ 0,-1 \}$
as soon as $d_c \ge 0$ \cite{OkonekT}.



\section{Dimensional reduction of DT theory to 4d Vafa-Witten theory}\label{sec:DTtoVW}

\subsection{Dimensional reduction of DT theory on an elliptic fibration}

Suppose that the 3-fold $X$ is an elliptic fibration over
a surface $S$. In particular, one can consider a trivial
fibration, $X = E \times S$. Then, we can interpret the
partition function of the Donaldson-Thomas theory on $X$
as counting BPS states in the T-dual system. Specifically,
under a T-duality along $E$, the D6-branes on $X$ turn into
D4-branes on $S$, the D0-brane charge becomes the charge
of the D2-brane on $E$, and so on.

For general gauge group $G$, D-brane charges can be viewed
as homology classes which are Poicar\'e duals of the obstruction classes
$o_i\in H^i(X,\pi_1(G))$. We consider compact connected Lie groups,
so $\pi_0(G)$ is trivial and the first obstruction class $o_1=0$.
The second obstruction class $o_2 \in H^2(X,\pi_1(G))$ is the
magnetic flux $m$. When $G=U(N)$, $m=c_1=ch_1 \in H^2(X,\Z)$ is the first Chern class
(which is also the first Chern character);
when $G=SO(N)\ (N>2)$, $m=w_2 \in H^2(X,\Z/2\Z)$, is the
second Steifel-Whitney class. In the rest of Section \ref{sec:DTtoVW},
we will consider $G=U(N)$, and the D-brane charges
are described by the Mukai vector as recalled in Section \ref{sec:reviewDT}.

In general, the resulting
configuration of D-branes is similar to the one we started with.
However, there are several cases when this duality can be
very useful and can relate Donaldson-Thomas invariants to
some invariants in the four-dimensional topological gauge theory.
The four-dimensional gauge theory is a twisted gauge theory
on $N$ D4-branes wrapped around the base $S$, obtained from
$N$ D6-branes in the original setup. This gauge theory is
a topological twist of $\CN=4$ super-Yang-Mills studied by Vafa and Witten.
It has the following spectrum \cite{VW}:

\medskip
\centerline{\vbox{\halign{\quad # & \quad # & \quad # \cr
bosonic: & 1-form $A$ & \cr
& self-dual 2-form $B_+$ & 3 \cr
& scalar $C$ & 1 \cr
& complex scalar $\phi$, $\bar \phi$ & 2 \cr
fermionic: & self-dual 2-forms $\chi$, $\tilde \psi$ & 3+3 \cr
& vectors $\psi$, $\tilde \chi$ & 4 + 4 \cr
& scalars $\eta$, $\zeta$ & 1+1 \cr }}}

\medskip
\noindent
where in the last column we again summarized the number of real
components.\footnote{Notice, that the total number of fermions is 16.}
Under a favorable set of
conditions\footnote{More specifically, when vanishing theorems hold.},
the partition function in this theory localizes on anti-self-dual
gauge fields (instantons):
\begin{equation} \label{eqn:ASD}
F_A^+ = 0
\end{equation}
and has the form similar to \eqref{ZDT}.

In general, the partition function of this four-dimensional
topological gauge theory can be interpreted as counting
BPS states with D2, D4, and D6 brane charges on $X$,
where all of these branes are wrapped on the elliptic fiber $E$.
Specifically, the dependence on the D2 and D4 brane charges
is through the gauge coupling $q = \exp( 2 \pi i \tau)$
and the 't Hooft magnetic fluxes $m$, while D6-brane charge
corresponds to the rank of the gauge group. In our case,
$G=U(N)$ and D2, D4, D6 brane charges are cycles dual
to $ch_2, ch_1, ch_0$, repsectively, and
the partition function of the Vafa-Witten theory can be interpreted
as a partition function of the Donaldson-Thomas theory,
\begin{equation}\label{eqn:ZDTVW}
Z_{DT}(X)_{ch_3=0} = Z_{VW}(S)= \sum_{k,m}\chi(\CM_{m,k}) t^m q^k
\end{equation}
where $\CM_{m,k}$ is the moduli space of solutions to
\eqref{eqn:ASD} with
$$
k = \int_S ch_2 (\CE), \quad \quad
m = ch_1 (\CE),\quad \quad N=ch_0(\CE).
$$
Notice that we have $c_1(T_X)=c_1(T_S)$ so the index \eqref{eqn:index} becomes
$$
\int_X c_1(T_S)\left(ch_0 ch_2 - \frac{ch_1^2}{2}+\frac{ch_0^2-1}{24}c_2(T_S)\right)
$$
which is zero when $ch_0, ch_1,ch_2$ are classes from $H^*(S,\Z)$, so this part
of the theory is balanced (in the sense of \cite{DM}). This is consistent with the fact that the
Vafa-Witten theory is balanced, even when $c_1(T_S)\neq 0$.

The partition function $Z_{DT}(X)=Z_{DT}(E\times S)$ is defined mathematically when
$S$ is a K3 surfaces. In this case, let
$$
\tN_{ch_1,ch_2,ch_3} = \int_{[\CM_{ch_1,ch_2,ch_3}]^\vir} 1
$$
be the Donaldson-Thomas (holomorphic Casson) invariants.
R. Thomas verified that in some cases $\tN_{m,k,0}= N^2 \chi(\CM_{m,k})$,
where the factor $N^2$ comes from the $N$-th root of unity in the Jacobian of $E$ \cite{Th}.


\subsection{Abelian Vafa-Witten theory}
In the abelian gauge theory, we have $ch_0(\CE)=1$.
In order for $ch_2 (\CE)$ to be non-trivial,
we need to enlarge the space of gauge connections to
include singular gauge fields (point-like instantons)
which mathematically can be described as ideal sheaves.
The partition function of the abelian gauge theory on $S$ is
\begin{equation}
Z(q,t) = {\vartheta_{\Gamma} (t,\tau,\bar \tau) \over \eta (\tau)^{\chi(S)} }
\label{zabelian}
\end{equation}
where the factor $\eta^{- \chi(S)}$ is the contribution from
the point-like instantons ($m=0$) and the theta-function comes from the lattice sum over
magnetic fluxes $m \in \Gamma$.

The moduli space of point-like instantons can be identified with Hilbert
scheme of points on $S$. According to Hirzebruch-Hofer \cite{HHofer},
the orbifold Euler number of the symmetric product $\Sym^k S$ is equal to
the Euler number of the Hilbert scheme $\Hilb^k S$ of $k$ points
in $S$, which is all we need.

G\"ottsche \cite{Gottsche} has shown:
$$
\sum_{k=0}^{\infty} q^k P_t (\Hilb^k S) =
\prod_{m=1}^{\infty} \prod_{i=1}^4 (1-(-t)^{2m-2+1} q^m)^{(-1)^{i+1} b_i(S)}
$$
In particular,
$$
\sum_{k \ge 0} \chi (\Hilb^k S) q^{k-\chi(S)/24} = \eta (\tau)^{-\chi(S)}
$$
or equivalently,
\begin{equation}\label{eqn:chiHilb}
\sum_{k \ge 0} \chi (\Hilb^k S) q^k =\prod_{n>0}\frac{1}{(1-q^n)^{\chi(M)}}
= q^{\frac{\chi(S)}{24}}\eta (\tau)^{-\chi(S)}
\end{equation}
This formula can be also obtained as a special case of the elliptic
genus for the symmetric product CFT with target space $\Sym^k S$ \cite{DMVVsym}.

\subsection{Examples and applications}

Below we consider two groups of examples:
$1)$ abelian gauge theory on $S$, and
$2)$ non-abelian $U(N)$ gauge theory with $m=0$.
The former is more relevant to the DT/GW correspondence,
while the latter gives us partition function of the DT theory
on elliptic fibration with $ch_1 (\CE) = 0$ and $ch_3 (\CE) = 0$.

\subsection{Abelian DT theory and GW theory on a trivial elliptic fibration}
Let $X=E\times S \to S$ be a trivial elliptic fiberation over
a smooth projective surface $S$, where $E$ is an elliptic curve.
In this subsection, we follow Edidin and Qin's computations of
DT and GW invariants of $X=E\times S$ in \cite{EQ}.

\subsubsection{Abelian DT theory}
Following \cite{MNOP1,MNOP2},
given $n\in \Z$ and $\beta\in H_2(X;\Z)$, let $I_n(X,\beta)$
be the moduli space of ideal sheaves $I_Z$ of 1-dimensional
closed subschemes $Z$ of $X$ satisfying
$\chi(\CO_Z)=n$ and $[Z]=\beta$.
The virtual dimension of $I_n(X,\beta)$ is $-\int_\beta K_S$.
Let $\beta_0=[E]\in H_2(X;\Z)$ be the fiber class.
The virtual dimension of $I_n(X,k\beta_0)$ is zero. Define
$$
\tN_{n,k}=\int_{[I_n(X,k\beta_0)]^{\vir}}1.
$$
The integral can be calculated by localization. The elliptic
curve $E$ is a group (compact torus). Let $E$ act on itself by
group multiplication and act on $S$ trivially. This gives a
torus action on $X=E\times S$ and torus actions on moduli
spaces $I_n(X,k\beta_0)$. The action has no fixed point unless
$n=0$. So
$$
\tN_{n,k}=0, \quad n\neq 0.
$$

When $n=0$, the moduli space $I_0(X,k\beta_0)$ can be identified with
$\Hilb^k S$, Hilbert scheme of $k$ points in $S$, which is
a smooth variety of complex dimension $2k$. The virtual
dimension of $I_0(X, k[E])$ is zero. Edidin and Qin proved
the following dimensional reduction
of Donaldson-Thomas theory on $X$ to Vafa-Witten theory on
$S$ in algebraic geometry:
\begin{fact} \cite[Lemma 3.4 (i)]{EQ}
The obstruction bundle over $I_0(X,d[E])\cong \Hilb^k S$
is isomorphic to the tangent bundle of $\Hilb^k S$.
\end{fact}
As a consquence \cite[Lemma 3.4 (ii)]{EQ},
$$
\tN_{0,k}=\int_{[I_0(X, k\beta_0) ]^\vir} 1
=\int_{\Hilb^k S} e(T_{\Hilb^k S})=\chi(\Hilb^k S).
$$
We have
$$
\sum_{k=0}^\infty \tN_{0,k}v^k=
\sum_{k=0}^\infty \chi(\Hilb^k S) v^k =\prod_{m>0}\frac{1}{(1-v^m)^{\chi(S)}}.
$$

\subsubsection{GW theory}\label{sec:S-E-GW}
For any $g,k\geq 0$ the moduli space
$\Mbar_{g,0}(X,k\beta_0)$ of genus $g$, degree $k\beta_0$ stable maps
to $X$ has virtual dimension zero. Define
$$
N_{g,k}=\int_{[\Mbar_{g,0}(X,k\beta_0)]^\vir}1.
$$
The moduli space $\Mbar_{g,0}(X,k\beta_0)$ can be identified with
$\Mbar_{g,0}(E,k[E])\times S$. The virtual dimension
of $\Mbar_{g,0}(X,k\beta_0)$ is zero, and the virtual
dimension of $\Mbar_{g,0}(E,k[E])$ is $2g-2$.
Let
$$
\pi_1:\Mbar_{g,0}(E,k[E])\times S\to \Mbar_{g,0}(E,k[E]),\quad
\pi_2:\Mbar_{g,0}(E,k[E])\times S\to S
$$
be projections to the first and second factors, respectively.
Let $\E\to \Mbar_{g,0}(E,k[E])$ be the Hodge bundle, and let
$\E^\vee$ denote its dual.
Then
$$
N_{g,k}=\int_{[\Mbar_{g,0}(X,k\beta_0)]^\vir} 1
= \int_{[\Mbar_{g,0}(E,k[E])]^\vir \times[S]}e(\pi_1^*\E^\vee\otimes \pi_2^* T_S)
$$
Again, $E$ acts on moduli spaces $\Mbar_{g,0}(X,k\beta_0)$, and
the action has no fixed points unless $g=1$. In this case, we have
$$
N_{1,k}=\int_{[\Mbar_{1,0}(E,k[E])]^\vir}1 \cdot \int_S e(T_S)=\chi(S) N_{1,k}^E
$$
where
$$
N_{1,k}^E=\int_{[\Mbar_{1,0}(E,k)]^{\vir}}1
$$
We have
$$
\sum_{k=1}^\infty N_{1,k}v^k = 
\chi(S) \sum_{k=1}^\infty N_{1,k}^E v^k.
$$
By \cite{BCOV,Dij},
$$
\exp\Bigl(\sum_{k=1}^\infty N^E_{1,k} v^k\Bigr)=\frac{1}{\prod_{m>0}(1-v^m)}.
$$
so
$$
\exp\Bigl(\sum_{k=1}^\infty N_{1,k}v^k\Bigr)
=\prod_{m>0}\frac{1}{(1-v^m)^{\chi(S)}}.
$$

\subsubsection{Non-compact surfaces}
The identity
\begin{equation}\label{eqn:chiS}
\sum_{k=0}^\infty \tN_{0,k} v^k
=\exp\Bigl(\sum_{k=1}^\infty N_{1,k}v^k\Bigr)=\prod_{m>0}\frac{1}{(1-v^m)^{\chi(S)}}.
\end{equation}
is also valid for certain noncompact surfaces. For example,
suppose that $S$ admits a $\C^*$ action such that the fixed point set
$S^{\C^*}$ is compact. Then one can use the torus action to define
residue DT invariants $\tN_{0,k}$ and GW invariants $N_{1,k}$.
We have
$$
N_{1,k}=\int_{[\Mbar_{1,0}(E,k[E])]^\vir}1 \cdot \int_{S^{\C^*}} \frac{e_{\C^*}(T_S)}{e_{\C^*}(N_{S^T/S})}
= N_{1,k}^E \chi(S^{\C^*})=N_{1,k}^E\chi(S).
$$
For example, $S$ can be the total space of a line bundle over a curve, or a resolution
of ADE singularities.

\subsubsection{Trivial elliptic fibration over a K3 surface}
Let $X=S\times E$, where $S$ be a projective K3 surface. Then any effective curve class 
$\beta\in H_2(X,\Z)$ is of the form
$$
\beta = \beta_1 + k\beta_0
$$
where $\beta_0=[E]$ as before, and $\beta_1$ is in the image of the
injective map $H_2(S,\Z)\to H_2(X,\Z)$ induced by the inclusion 
$$
S = S\times \{p_0\}\hookrightarrow S\times E = X
$$
where $p_0$ is some point in $E$.

\bigskip

\noindent
$\underline{\beta_1=0}$: by Section \ref{sec:S-E-GW}, $N_{g,k\beta_0}^{S\times E}=0$ unless $g=1$, and
$$
\exp\left(\sum_{k=0}^\infty N^{S\times E}_{g=1, k\beta_0} v^k\right) =\prod_{m>0}\frac{1}{(1-v^m)^{24}}.
$$

\noindent
$\underline{\beta_1>0}$ (viewed as a class in $H_2(S,\Z)$): we have (see e.g. \cite{MP,ObP})
$$
[\Mbar_{g,0}(S\times E,\beta_1+k\beta_0)]^\vir =0,
$$
so
$$
N_{g,\beta_1+k\beta_0}^{S\times E}=0.
$$
Therefore, Gromov-Witten theory of $S\times E$ is almost trivial.
However, $S\times E$ has very interesting {\em reduced} Gromov-Witten theory
\cite{ObP, ObS}.


\subsection{Half K3}

\subsubsection{Construction}
Let $[\XX]$ be homogeneous coordinates on $\cP^2$. Let $F_0(\XX)$
and $F_1(\XX)$ be two homogeneous polynomials of degree 3 such that
\begin{eqnarray*}
D_0&=&\{ [\XX]\in\cP^2 \mid F_0(\XX)= 0\},\\
D_1&=&\{ [\XX]\in\cP^2 \mid F_1(\XX)=0\}
\end{eqnarray*}
are two smooth cubic curves intersecting transversally at 9 points
$q_1,\ldots, q_9$. $D_0$, $D_1$ are elliptic curves.

Consider a pencil of elliptic curves in $\cP^2$ containing $D_0$
and $D_1$:
$$
B_9=\{  ([\z], [\XX])\in \cP^1\times \cP^2\mid Z_0 F_0(\XX)+ Z_1
F_1(\XX)=0\}.
$$

Define $p_1:\cP^1\times \cP^2\to \cP^1$ and $p_2:\cP^1\times \cP^2\to \cP^2$
to be the projections to the first and second factors, respectively.
And, let $\pi_1:B_9\to \cP^1$ and $\pi_2:B_9\to \cP^2$ be the restrictions of these projections.
Then $\pi_1:B_9\to \cP^1$ is an elliptic fibration, and $\pi_2: B_9\to
\cP^2$ is blowup of $\cP^2$ at $q_1,\ldots,q_9$.
The exceptional divisor $E_i=\pi_2^{-1}(q_i)$ is a section of $\pi_1: B_9 \to \cP^1$
and can be identified with $\cP^1\times \{q_i\}\subset \cP^1
\times \cP^2$.


\subsubsection{Basic topology}

Let $H_1, H_2\in H^2(\cP^1\times \cP^2;\Z)$ be the pull back of
the hyperplane class under $p_1,p_2$, respectively. Then
$$
H^*(\cP^1\times \cP^2;\Z)=\Z[H_1,H_2]/ (H_1^2, H_2^3).
$$
Let $h_1, h_2\in H^2(B_9;\Z)$ be the restriction of $H_1, H_2$
respectively. Then $h_1$ is the Poincar\'{e} dual of the homology
class of a fiber $F$ of $\pi_1:S\to \cP^1$, and $h_2$ is the
Poincar\'{e} dual of the homology class of a curve $E_0\subset S$ such
that $\pi_1|_{E_0}: E_0\to \cP^1$ is a degree 3 cover. We have
$E_0\cdot E_0=1$ and the genus of $E_0$ is zero. By adjunction
formula, we have $c_1(T_{B_9})=h_1$. So the pair $(B_9,F)$ is relative
Calabi-Yau in the sense that $K_{B_9} + F=0$. In particular,
$h^{2,0}(B_9)=h^{0,2}(B_9)=0$.

We have
$$
\chi(B_9)=\int_{B_9}c_2(T_{B_9})=\int_{\cP^1\times \cP^2}
\frac{(1+H_1)^2(1+H_2)^3}{1+H_1+3H_2}(H_1+3H_2) =12.
$$
Alternatively,
$$
\chi(B_9)=\chi(\cP^2)-9\chi(\mathrm{point})+ 9\chi(\cP^1) =3-9+18=12
$$

By the Lefshetz hyperplane theorem, $H_1(B_9;\Z)=0$. So the Hodge
diamond of $B_9$ is given by:
$$
\begin{array}{ccccccc}
   &   &     & 1 &    &    &   \\
   &   &  0 &    & 0 &    &   \\
   & 0&     & 10 &    & 0 &   \\
   &    & 0 &    & 0 &  & \\
   &    &    & 1 &    &   &
\end{array}
$$

For $i=0,\ldots,9$, let $e_i\in H_2(B_9;\Z)$ be the homology class
represented by the curve $E_i$, respectively. Then
$$
H_2(B_9;\Z)=\bigoplus_{i=0}^9 \Z e_i,
$$
and the intersection form is given by
$$
e_i\cdot e_j = g_{ij},\ \  g_{ij}=\mathrm{diag}(1,-1,\ldots, -1).
$$
so that
$$
H_2 (B_9,\Z) = \Gamma^{1,9} 
$$
is the unique 10 dimensional odd unimodular lattice of signature $(1,9)$.

In terms of $e_i$, the fiber class $[F]$ and the base class $[B]$
have the form
$$
[F] = 3 e_0 -\sum_{i=1}^9 e_i, \quad\quad [B] = e_9
$$
such that
$$
[F] \cdot [B] = 1 \quad\quad [B] \cdot [B] = -1 \quad\quad
[F]\cdot [F] = 0
$$
We have an orthognal decomposition
$$
\Gamma^{1,9}=\Gamma^{1,1}\oplus E_8(-1)
$$
where 
\begin{itemize}
\item $\Gamma^{1,1}$ is the sublattice generated by $[F]$ and $[B]$, 
which is the unique 2 dimensional odd unimodular lattice of signature
$(1,1)$, and
\item $E_8(-1)$ is the sublattice generated by $-e_0+e_1+e_2+e_3,
e_1-e_2, e_2-e_3,\ldots, e_7-e_8$,
which is the negative $E_8$ lattice.
\end{itemize}
\subsubsection{Local GW invariants from topological YM theory}

Recall that $K_{B_9}=-F$.
Suppose that
$$
\beta\in H_2(X_{B_9};\Z)\cong H_2(B_9;\Z)
$$
and $\beta\cdot[F]>0$. Then
$$
\Mbar_{g,0}(X_{B_9},\beta)=\Mbar_{g,0}(B_9,\beta).
$$

The virtual dimension of $\Mbar_{g,0}(X_{B_9};\beta)$ is zero, while
the virtual dimension of $\Mbar_{g,0}(B_9;\beta)$ is
$\beta\cdot[F]+g-1$. Define
\begin{equation}\label{eqn:Bnine}
N^{X_{B_9}}_{g,\beta}=\int_{[\Mbar_{g,0}(X_{B_9},\beta)]^{\vir}} 1=
\int_{[\Mbar_{g,0}(B_9,\beta)]^{\vir}}e(V_{g,\beta})
\end{equation}
where $V_{g,\beta}$ is a vector bundle over
$\Mbar_{g,0}(B_9;\beta)$ whose fiber over $[f:C\to B_9]$ is
$$
H^1(C,f^*\CO_{B_9}(-F)).
$$
Note that $H^0(C,f^*\CO_S(-F))=0$ since $[C]\cdot (-[F])<0$. By
Riemann-Roch, the rank of $V_{g,\beta}$ is $\beta\cdot[F]+g-1$.

Introduce the Gromov-Witten potential,
$$
\CF(\lambda;t) := \sum_g\lambda^{2g-2} \CF_g(t),\quad \quad
\CF_g(t)=\sum_{\b}N_{g,\b}^{X_{B_9}}  e^{2 \pi i t \cdot \b}
$$
where $t = (\phi, \tau, \vec m)$ are the K\"{a}hler moduli, such that
$\phi$ and $\tau$ correspond to the base and fiber moduli, and
$m_i$, $i=1,\ldots,8$ denote the remaining $E_8$ moduli.

Since $B_9$ is elliptic, one can use the fiberwise T-duality to
relate the counting of holomorphic curves on $B_9$ to the counting
of instantons on $B_9$. As a result, one finds a relation
\begin{equation}
\label{fviazn} \CF_0(t) = \sum_{n=1}^{\infty} q^{n/2} Z_n (m_i ; \tau)
e^{2\pi i n \phi}
\end{equation}
where $q = e^{2\pi i \tau}$ and $Z_n (m_i ; \tau)$ is the
partition function of the topological $\CN=4$ $U(n)$ Yang-Mills on
$B_9$. It is defined as a generating function
$$
Z_n := \sum_k \chi (\CM_k) q^k
$$
where $\CM_k$ is the moduli space of instantons on $S$ with
instanton number $k$. Thus, \eqref{fviazn} relates the counting of
holomorphic curves in the class $\beta = n[B] + k[F]$ to the
counting of $U(n)$ instantons with instanton number $k$. Let us
consider some explicit examples.

\noindent$\underline{n=0}:$ In this case, the rank of the gauge
group is zero and there are no Yang-Mills instantons. Therefore,
$$
N_{g,\b}^{X_{B_9}} = 0
$$
for $\b = k [F]$.

\noindent$\underline{n=1}:$ In the abelian Yang-Mills theory, the
only non-trivial instantons are point-like. The moduli
space of point-like instantons of charge $k$ can be identified with
$\Hilb^k B_9$, the Hilbert scheme of $k$ points on the surface $B_9$. This
suggests $Z_1 \sim \eta(q)^{-1/12}$.

In general,
$$
 q^{1/2}Z_1(m_i;\tau)=\frac{\theta_{E_8}(\tau;m_i)q^{1/2}}{\eta(q)^{12}}
=\theta_{E_8}(\tau;m_i)\prod_{m>0}\frac{1}{(1-q^m)^{12}}.
$$
The precise definition of $\theta_{E_8}(\tau;m_1,\ldots,m_8)$
can be found in \cite[section 2.6]{Ho}. It can be written as
$$
\theta_{E_8}(\tau;m_1,\ldots,m_8)=\sum_{k_i\in
\frac{1}{2}\Z} A_{k_1,\ldots,k_8}(q) y_1^{k_1}\cdots y_8^{k_8}
$$
where $y_i= e^{2\pi i m_i}$, $A_{0,\ldots,0}(q)= 1$. Hence, one finds
\begin{equation}\label{eqn:oneB}
\begin{aligned}
& \sum_{k=0}^\infty N_{0,[B]+k[F]}^{X_{B_9}} q^k =\prod_{m>0}\frac{1}{(1-q^m)^{12}} \\
=&1+12q+90 q^2+ 520 q^3 + 2535 q^4 + 10908 q^5
+\cdots
\end{aligned}
\end{equation}

In fact, the duality prediction \eqref{eqn:oneB} has been verified mathematically.
The rank of $V_{0,[B]+k[F]}$ in \eqref{eqn:Bnine} is zero, so
$$
N^{X_{B_9}}_{0,[B]+k[F]}=\int_{[\Mbar_{0,0}(B_9;[B]+k[F])]^\vir}1
=N^{B_9}_{0,[B]+k[F]}
$$
Let $N_{g,n}= N^{B_9}_{g,[B]+(g+n)[F]}$, where
$N^{B_9}_{g,\beta}$ is defined by \eqref{eqn:NS}.
Bryan and Leung \cite{BL} showed that
$$
\sum_{n=0}^\infty N_{g,n} q^n= \left(\sum_{k=1}^\infty
k(\sum_{d|k} d) q^{k-1} \right)^g \prod_{m>0}\frac{1}{(1-q^m)^{12}}.
$$
The case $g=0$ gives exactly the formula \eqref{eqn:oneB}.

\bigskip

\noindent$\underline{n \ge 2}:$ In general, $Z_n$ have the
following structure \cite{MNVW}:
\begin{equation}
Z_n (0;\tau) = \eta^{-12n} P_{6n-2} (E_2 (\tau), E_4 (\tau), E_6
(\tau))
\end{equation}
where $P_{6n-2}$ is a quasi-homogeneous polynomial of weight
$6n-2$, and $E_{2m}$ are the Eisenstein modular forms of weight $2m$.

In general $q^{n/2}Z_n(0;\tau)$ gives us
\begin{equation}
\sum_k \left(\sum_{\beta\cdot [B]=k, \beta\cdot[F]=n}N_{0,\beta}^{X_{B_9}}\right)q^k.
\end{equation}
instead of
\begin{equation}
\sum_{k=0}^\infty N_{0,n[B]+k[F]}^{X_{B_9}} q^k
\end{equation}
For example,
\begin{equation}
\begin{aligned}
q^{1/2}Z_1(0;\tau)&=\frac{\theta_{E_8}(\tau;0)q^{1/2}}{\eta(q)^{12}}
=\frac{E_4(\tau)}{\prod_{m>0}(1-q^m)}\\
&=1+12q+330q^2+3400 q^3
+26295q^4+161628 q^5 +\cdots
\end{aligned}
\end{equation}
which is different from \eqref{eqn:oneB}.

We have
\begin{equation}\label{eqn:twoB}
q Z_2(0;\tau)=
-\frac{1}{8}+\frac{18441}{2}q+ 673760 q^2+ \frac{82133595}{4} q^3 +\cdots
\end{equation}

Gromov-Witten invariants of $X_{B_9}$, i.e. local Gromov-Witten
invariants of the half K3, are studied in
\cite{HSTa}, \cite{Sa}, \cite{KMW}, etc.



\subsection{Dimensional reduction of DT theory on a local surface}

A closely related problem appears when we study the total space, $X_S$,
of the canonical line bundle of a complex algebraic surface $S$.
We will be mainly interested in enumerative invariants of $X_S$
when $S$ is a Fano surface.
In this case, any non-constant holomorphic map to $X_S$
factors through the zero section $S \to X_S$,
so the Gromov-Witten invariants of $X_S$ can be viewed
as local Gromov-Witten invariants of $S$ in a Calabi-Yau three-fold.
Similarly, Donaldson-Thomas gauge theory on $X_S$
effectively reduces to a four-dimensional gauge theory on $S$.

\subsubsection{The reduced deformation complex}
In order to see this more explicitly,
let us consider the deformation complex \eqref{sixdcplx}
for the total space $X_S$, of the canonical line bundle of $S$.
Let $v$ be a local complex coordinate on the fiber of $X_S$.
Then, a $(p,q)$-form on $X_S$ which is constant along
the fiber locally can be written as
$$
\omega^{(p,q)} = \tilde \omega^{(p,q)}
+ \tilde \omega^{(p,q-1)} \wedge d \bar v
+ \tilde \omega^{(p-1,q)} \wedge d v
+ \tilde \omega^{(p-1,q-1)} \wedge dv \wedge d \bar v
$$
where we use $\tilde \omega^{(p,q)}$ to denote a (pull-back)
of a $(p,q)$-form on $S$. Let $\hat{\Omega}^{0,k}$ be the
space of $(0,k)$-forms on $X_S$ which are constant along the fiber.
Then it decomposes
into the following spaces of anti-holomorphic forms on $S$:
\begin{equation}
\begin{aligned}
\hat{\Omega}^{0,3}&= \Omega^{0,2}_S \\
\hat{\Omega}^{0,k} &= \Omega^{0,k}_S \oplus \Omega^{0,k-1}_S \quad\quad k=1,2 \\
\hat{\Omega}^{0,0} &= \Omega^{0,0}_S
\end{aligned}
\end{equation}
Hence, in terms of $(0,k)$-forms on $S$, the complex \eqref{sixdcplx} reads
\begin{equation}
0 \to \Omega^{0,0} \to
\Omega^{0,0} \oplus \Omega^{0,1} \oplus \Omega^{0,2} \to
\Omega^{0,1} \oplus \Omega^{0,2} \to 0
\label{reducedcplx}
\end{equation}

We now describe the differentials in \eqref{reducedcplx}.
We assume that the fields are contant along the fiber
of $X_S\to S$. Then $A\in \Omega^1_{X_S}(\End\CE)$ can be decomposed as
\begin{equation}\label{eqn:A}
A= A_2 + C dv + \bar{C} d\bar{v},
\end{equation}
where $v$ is the local coordinate in the fiber of $X\to S$,
$$
A_2\in \Omega^1_S(\End\CE),\quad C\in \Omega^0_S(\End\CE\otimes K_S^{-1}),\quad
\bar{C}\in \Omega^0_S(\End\CE\otimes K_S).
$$
and $\bar{\varphi}$ can be written as
\begin{equation}\label{eqn:phi}
\bar{\varphi}=\bar{\varphi}_2\wedge d\bar{v}
\end{equation}
where $\bar{\varphi}_2\in \Omega^{0,2}_S(\End\CE\otimes K_S)$.

Equations \eqref{eqsprime} can be written as
\begin{equation}\label{eqn:FC}
F_{A_2}^{0,2}=[C,\bar{\varphi}_2]
\end{equation}
\begin{equation}\label{eqn:Cphi}
\dbar_{A_2} \bar{C} =\dbar_{A_2}^\dagger \bar{\varphi}_2
\end{equation}
\begin{equation}
\Tr_\omega F_{A_2}^{1,1} + [C,\bar{C}] + * [\varphi_2,\bar{\varphi}_2]=\ell I.
\end{equation}
where  $\omega$ is the K\"{a}hler form on $S$.

Let $V=\End\CE$.
The deformation complex at $(C,A_2,\varphi_2)$ is given by
$$
0\to \Omega_S^{0,0}(V)\stackrel{D_1}{\to}
\Omega_S^{0,0}(V\otimes K_S) \oplus\Omega_S^{0,1}(V)\oplus
\Omega_S^{0,2}(V\otimes K_S)
\stackrel{D_2}{\to}
 \Omega_S^{0,1}(V\otimes K_S)\oplus \Omega_S^{0,2}(V)\to 0
$$
where
\begin{eqnarray*}
D_1(\eta)&=&  ([\bar{C},\eta], \dbar_{A_2}\eta, [\bar{\varphi}_2,\eta])\\
D_2(\delta\bar{C},\delta A_2^{0,1},\delta \bar{\varphi}_2)&=&
\left([\bar{C},\delta A_2^{0,1}]+\dbar_{A_2}(\delta\bar{C})
-\dbar_{A_2}^\dagger (\delta\bar{\varphi}_2),
 \dbar_{A_2} (\delta A_2^{0,1}) - [C,\delta\bar{\varphi}_2]\right)
\end{eqnarray*}
This complex has to be identified with a deformation complex
in a suitable four-dimensional gauge theory.


\subsubsection{Three twists of $\CN=4$ super-Yang-Mills}

There are three different twists of $\CN=4$ super-Yang-Mills
in four dimensions, which can be conveniently characterized by the following homomorphisms \cite{VW}:

\begin{enumerate}

\item $SU(2)_L \times SU(2)_R \times SU(4)_I \to SU(2)_L'
\times SU(2)_R \times SU(2)_F$, under which ${\bf 4} \to ({\bf 2},
{\bf 1}) \oplus ({\bf 2}, {\bf 1})$

\item $SU(2)_L \times SU(2)_R \times SU(4)_I \to SU(2)_L'
\times SU(2)_R \times SU(2)_F \times U(1)$, so that ${\bf 4} \to
({\bf 2}, {\bf 1}) \oplus ({\bf 1}, {\bf 1}) \oplus ({\bf 1}, {\bf
1})$

\item $SU(2)_L \times SU(2)_R \times SU(4)_I \to SU(2)_L'
\times SU(2)_R' \times U(1)$, so that ${\bf 4} \to ({\bf 2}, {\bf
1}) \oplus ({\bf 1}, {\bf 2})$

\end{enumerate}

\noindent
and which relate the rotation symmetry and R-symmetry in the original theory and its twisted (topological) version.
The first twist leads to a TQFT known as the Vafa-Witten theory \cite{VW}.
The second is the Donaldson-Witten twist that leads to non-abelian monopole equations \cite{LabastidaL}.
And the last twist is sometimes called Marcus twist \cite{Marcus} or GL twist, for its connection to
the geometric Langlands program \cite{KW}.

These three TQFTs involve (though in a somewhat different way)
certain moduli spaces: $1)$ moduli space of self-dual gauge
connections; $2)$ moduli space of adjoint non-abelian monopoles;
$3)$ moduli space of complexified flat gauge connections.
The second theory based on non-abelian monopole equations
has only one topological supersymmetry, $\CN_T = 1$.
On the other hand, theories $1)$ and $3)$ have $\CN_T = 2$
topological supersymmetry.

For example, in the Vafa-Witten $\CN=4$ twisted gauge theory,
the deformation complex is \cite{VW}:
\begin{eqnarray}
0 \to \Omega^0_S( {\rm ad} P)
&\to& \Omega^0_S({\rm ad} P) \oplus \Omega^1_S({\rm ad} P) \oplus \Omega^{2,+}_S({\rm ad} P) \cr
&\to& \Omega^1_S({\rm ad} P) \oplus \Omega^{2,+}_S({\rm ad} P) \to 0
\label{vwcplx}
\end{eqnarray}

For comparison, let us write the deformation complex in the generalization
of $\CN=2$ Seiberg-Witten theory, where the monopole fields take
values in $\Gamma (S, S^+ \otimes E)$,
$$
0 \to \Omega^0 (g_E) \to \Omega^1 (g_E) \oplus \Gamma (S, S^+ \otimes E)
\to \Omega^{2,+} (g_E) \oplus \Gamma (S, S^- \otimes E) \to 0
$$
Here, $g_E$ denotes the representation of the Lie algebra $g = Lie (G)$
associated to a representation $R$. If $P \to S$ is the principal
$G$-bundle over $S$, and $V$ is a vector space associated with
representation $R$, we can form the associated vector bundle $E = P \times_G V$.

In order to compare the deformation complex of the six-dimensional
gauge theory on $X_S$ with deformation complexes in familiar
four-dimensional topological gauge theories $(1)$-$(3)$, we now
assume that $S$ is K\"{a}hler, and use the famous result of Donaldson.
Namely, for a bundle $E$ over a Kahler surface $S$, we can identify
the moduli space of irreducible ASD connections with the set
of equivalence classes of stable holomorphic bundles $\CE$
which are topologically equivalent to $E$. We fix the holomorphic
line bundle $\det\CE$ so that $G=SU(n)$.
There are natural isomorphisms
\begin{eqnarray}
\Omega^1 (g_E) & \cong & \Omega^{0,1} (\End_0\CE) \cr
\Omega^0 (g_E) \oplus \Omega^{2,+} (g_E) & \cong &
\Omega^{0} (\End_0\CE) \oplus \Omega^{0,2} (\End_0\CE)
\end{eqnarray}
where $\End_0 \CE$ denotes the trace-free endomorphisms of
the stable holomorphic vector bundle $\CE$.
Using such isomorphisms, we can write the deformation
complex \eqref{reducedcplx} as
$$
0 \to \Omega^0 \to \Omega^0 \oplus \Omega^1 \oplus \Omega^{2,+} \to
\Omega^1 \oplus \Omega^{2,+} \to 0
$$
where elements of all these spaces are valued in $g_E$
associated with $\End_0 (\CE)$.
Now, this complex does look like the deformation complex \eqref{vwcplx}
in Vafa-Witten theory!

On a K\"{a}hler surface $S$, we further have isomorphisms
$$
S^+ = \Omega^{0,0} \oplus \Omega^{0,2}
$$
$$
S^- = \Omega^{0,1}
$$
$$
\Omega_{\mathbb{C}}^{2,+}  = \Omega^{2,0}
\oplus \Omega^0  \cdot \omega \oplus \Omega^{0,2}
$$
As pointed out in \cite{Th}, there is a simple way to relate
the Vafa-Witten theory on $S$ and the Donaldson-Thomas theory on
$X_S$. Given a stable bundle $\CE$ on $S$, let $i_*\CE$ be the
pushforward torsion sheaf on $X_S$ supported on $S$. This gives
an inclusion $i_*: \CM(S) \longrightarrow \CM(X_S)$, whose image
is a connected component of $\CM(X_S)$. The restriction of
the deformation theory on $\CM(X_S)$ to $i_*\CM(S)$ coincides
with that on $\CM(S)$.


\section{Dimensional reduction of DT theory to 2d $\sigma$-model} \label{sec:DTtoSigma}

\subsection{Topological reduction}
The second duality we consider is similar to the topological
reduction of four-dimensional gauge theory studied in \cite{BJSV}.
Namely, we consider $X$ to be of the form $X = \Sigma_\ell \times S$,
where $\Sigma_\ell$ is a Riemann surface of genus $\ell$
and $S$ is a K\"{a}hler surface. Let $g_\Sigma$ and
$g_S$ be the metrics of $\Sigma$ and of $S$, respectively.
We first consider a general rescaling
$$
g_\Sigma \to t^p g_\Sigma,\quad
g_S\to t^q g_S
$$
where $t$ is a positive number which tends to zero. Then the Lagrangian
density of the Yang-Mills action scales as
$$
F_X \wedge * F_X \to  t^p F_S \wedge * F_S
+ 2 t^q(d_S A_\Sigma -D_\Sigma A_S)\wedge *(d_S A_\Sigma - D_\Sigma A_S)
+ t^{-p+2q} F_\Sigma\wedge * F_\Sigma.
$$

When $(p,q)=(0,1)$, we have $g_\Sigma\to g_\Sigma, \quad g_S \to t g_S$ ($S$ is ``small''), and
$$
F_X \wedge * F_X \to  F_S \wedge * F_S
+ 2t (d_S A_\Sigma -D_\Sigma A_S)\wedge *(d_S A_\Sigma - D_\Sigma A_S)
+ t^2 F_\Sigma\wedge * F_\Sigma.
$$
The resulting theory is the 4d Vafa-Witten theory discussed in section \ref{sec:DTtoVW}.

This section concerns the case $(p,q)=(-1,0)$. In this case,
we have
$$
g_\Sigma \to \frac{1}{t}g_\Sigma,\quad g_S\to g_S \quad \textup{($\Sigma$ is ``large'') },
$$
$$
F_X \wedge * F_X \to  \frac{1}{t} F_S \wedge * F_S
+ 2 (d_S A_\Sigma -D_\Sigma A_S)\wedge *(d_S A_\Sigma - D_\Sigma A_S)
+ t F_\Sigma\wedge * F_\Sigma.
$$
As in \cite{BJSV}, the resulting theory on $\Sigma_{\ell}$
is $\CN=2$ topological sigma-model, whose target space is
the moduli space of solutions to the equations in
the {\it four-dimensional} gauge theory on $S$.
As we explained in Section \ref{sec:DTtoVW}, these are
anti-self-dual gauge fields on $S$ (when vanishing theorems
hold).
Therefore, after the topological reduction we obtain
$\CN=2$ sigma-model on $\Sigma_{\ell}$ with
the target space $\CM_{ASD} (S)$.

Since the six-dimensional theory
is topological, it should not matter whether we reduce on
$S$ or on $\Sigma_{\ell}$. This means that the partition
function $Z_{DT} (X)$ should be equal to the partition
function of the four-dimensional gauge theory on $S$,
as well as to the partition function of the $\CN=2$ sigma-model
on $\Sigma_{\ell}$,
$$
Z_{DT} (X) = Z_{4D} (S) = Z_{2D} (\Sigma_{\ell})
$$

\subsection{Moduli spaces of instantons and correlation functions}

We consider a general gauge group $G$ which can be any compact connected Lie group.
Let $\CM_{\mathrm{ASD}}(S,G)$ denote the moduli space of gauge equivalence classes of
ASD $G$-connections on $S$, which are solutions to \eqref{eqn:ASD}.
Let  $\CM_{\mathrm{HYM}}(\Sigma_\ell\times S, G)$ denote the moduli space
of gauge equivalence classes of HYM $G$-connections on
$\Sigma_\ell\times S$, which are solutions  to \eqref{eqn:HYM}.
Each of the spaces $\CM_{ASD}(S,G)$ and $\CM_{HYM}(\Sigma_\ell\times S, G)$
can be represented as a union of components labelled by the topological type of the
underlying topological principal $G$-bundle.

When $G$ is semi-simple, the ASD equation \eqref{eqn:ASD} is equivalent to
the 2-dimensional version of the HYM equations \eqref{eqn:HYM}.
When $G$ is not semi-simple, there is
a topological obstruction to the existence of ASD connections.
Let $H$ be the connected component of $Z(G)$, the center of $G$. Then
$H$ is a compact torus. Let $G_{ss}=[G,G]$ be the commutator subgroup,
which is also the maximal connected semismiple subgroup of $G$. There
is a finite cover $H\times G_{ss} \to G$ given by $\quad (h,g)\mapsto hg$.
This is a group homomorphism. The kernel is isomorphic to
$D=H\cap G_{ss}$ which is a finite subgroup of $H$. Note that $G_{ss}$ is
a normal subgroup of $G$. Consider the exact sequence of abelian groups
$$
1 \to \pi_1(G_{ss})\to \pi_1(G)\to \pi_1(G/G_{ss})\to 1
$$
where $G/G_{ss}= H/D \cong \Z^{\oplus \dim_\R H}$. We have
$$
\cdots \to H^2(S,\pi_1(G_{ss}) ) \to  H^2(S,\pi_1(G)) \stackrel{\deg}
\to H^2(S,\Z)^{\oplus \dim_\R H}  \to \cdots
$$
Let $m=o_2(P)\in H^2(S,\pi_1(G))$ be the magnetic flux of
a principal $G$-bundle $P$, so that $\deg(m)=(m_1,\ldots,m_r)$, where
$m_i \in H^2(S;\Z)$ and $r=\dim_\R H$. Then $P$ admits an ASD
connection only if $m_i \cdot [k_0]=0$, where $k_0$ is the K\"{a}hler
class. Given $c \in H^2(S;\Z)\setminus\{0 \}$, let
$$
\Gamma_c=\{  k \in H^{1,1}(S)\mid c\wedge k\wedge k=0\}.
$$
Then $\Gamma_c$ is a real codimension-1 subspace of $H^{1,1}(S)$.
Choose $k_0$ such that
$$
[k_0]\notin \bigcup_{c\in H^2(S;\Z), c\neq 0} \Gamma_c
$$
Then $m_i \cdot [k_0]=0$ implies $m_i=0$.
In this case $P$ admits an ASD connection only if
$$
\deg(P)\stackrel{\mathrm{def}}{=}\deg(o_2(P))=0\in H^2(S;\Z)^{\oplus r}.
$$
For example, when $G=U(N)$, we have $G_{ss}=SU(N)$ and
$o_2(P)=\deg(P)=c_1(P)$. When $G=SO(N)\ (N>2)$, we have $G_{ss}=SO(N)$,
$o_2(P)=w_2(P)$, and $\deg(P)=0$.

The complexification $G^\C$ of $G$ is a connected reductive algebraic group
over $\C$. The moduli space $\CM_{\mathrm{ASD}}(S,G)$ can be identified
with the moduli space $\CM(S,G^\C)$ of $S$-equivalence classes of
semi-stable holomorphic principal $G^\C$-bundles of degree zero on $S$ \cite{Do},
and the moduli space
$\CM_{\mathrm{HYM}}(\Sigma_\ell\times S,G)$ can be identified with the moduli space
$\CM(\Sigma_\ell\times S,G^\C)$
of $S$-equivalence classes of semi-stable holomorphic principal $G^\C$-bundles  over
$\Sigma_\ell\times S$ \cite{UY, AB}.

The partition function of the $\CN=2$ sigma-model on $\Sigma_\ell$ localizes
on $\Mor(\Sigma_\ell,\CM(S,G^\C))$, the space of holomorphic maps from
$\Sigma_\ell$ to $\CM(S,G^\C)$. The partition function of the Donaldson-Thomas
theory on $\Sigma_\ell \times S$ localizes on $\CM(\Sigma_\ell\times S,G^\C)$.
We now construct a map between the two moduli spaces
$\Mor(\Sigma_\ell,\CM(S,G^\C))$ and $\CM(\Sigma_\ell\times S,G^\C)$.
Let $\CP \to \CM(S,G^\C)\times S$ be the universal principal $G^\C$-bundle.
(It does not exist as a usual principal bundle, see e.g. \cite[Section 7.1]{KW}.)
Given a holomorphic map $u:\Sigma_\ell\to \CM(S,G^\C)$, we get a holomorphic
principal $G^\C$-bundle $(u\times \mathrm{id}_S)^*\CP$ over $\Sigma_\ell\times S$.
We get a map
\begin{equation}\label{eqn:PhiMM}
\Phi:\Mor(\Sigma_\ell,\CM(S,G^\C)) \longrightarrow \CM(\Sigma_\ell\times S,G^\C),\quad
u\mapsto (u\times \mathrm{id}_S)^*\CP
\end{equation}
which is birational.

Let $\CO_1,\ldots,\CO_k$ be observables in the $\sigma$-model on $\Sigma_\ell$.
Our topological reduction implies that there are corresponding observables
$\tilde{\CO}_1,\ldots,\tilde{\CO}_k$ in the DT theory on $\Sigma_\ell\times S$ such
that
\begin{equation}\label{eqn:OtO}
\langle \CO_1 \cdots \CO_k \rangle
= \langle \tilde{\CO}_1 \cdots \tilde{\CO}_k \rangle
\end{equation}
Geometrically, these correlation functions should be interpreted as intersection
numbers on moduli spaces $\Mor(\Sigma_\ell,\CM(S,G^\C))$ and
$\CM(\Sigma_\ell\times S,G^\C)$. To define these intersection numbers, we need
to compactify these moduli spaces.

When $S$ is a projective surface, the moduli spaces $\CM(S,G^\C)$ and $\CM(\Sigma_\ell\times S,G^\C)$
can be compactified  using algebraic geometry. When $G=U(N)$, $G^\C=GL(N,\C)$,
these moduli spaces can be identified with moduli spaces of isomorphism classes of
polystable holomorphic vector bundles, and the compactified moduli spaces are obtained by
including semistable torsion free sheaves.  For general $G$, the compactified moduli spaces
$\overline{\CM}(S,G^\C)$ and $\overline{\CM}(\Sigma_\ell\times S,G^\C)$ are obtained by including
semistable {\em principal sheaves} \cite{GS} or {\em singular principal bundles} \cite{Sch}.
These compactified moduli spaces are projective but may have singularities.

To compactify $\Mor(\Sigma_\ell,\CM(S,G^\C)$, we first compactify the target
and get a partial compactification $\Mor(\Sigma_\ell,\Mbar(S,G^\C))$
which consists of morphisms from $\Sigma_\ell$ to $\Mbar(S,G^\C)$.
The final compactification $\overline{\Mor}(\Sigma_\ell,\Mbar(S,G^\C))$
is obtained by including stable maps from a curve with a root component
isomorphic to $\Sigma_\ell$ and bubble components which are spheres.
Ideally, \eqref{eqn:PhiMM} extends to a morphism
$\Phi:\overline{\Mor}(\Sigma_\ell,\CM(S,G^\C)) \longrightarrow \Mbar(\Sigma_\ell\times S,G^\C)$.
The inverse is a rational map: a semi-stable $G^\C$-bundle $P\to \Sigma_\ell\times S$ representing
an element in $\CM(\Sigma_\ell\times S, G^\C)$ defines a morphism $\Sigma_\ell \to \CM(S, G^\C)$.
$$
\tilde{\CO}_i\in H^*(\Mbar(\Sigma_\ell\times S,G^\C)),\quad
\CO_i=\Phi^*\tilde{\CO}_i \in  H^*(\overline{\Mor}(\Sigma_\ell,\Mbar(S,G^\C)).
$$
The correlation functions \eqref{eqn:OtO} can be interpreted
as intersection numbers of the cohomology classes
$\CO_i$'s and $\tilde{\CO}_i$'s. For a singular moduli space
one needs a virtual fundamental class to define intersection
theory. The topological reduction suggests that there exist
virtual fundamental classes such that
$$
\Phi_*[\overline{\Mor}(\Sigma_\ell,\Mbar(S,G^\C)]^\vir =
[\Mbar(\Sigma_\ell\times S,G^\C)]^\vir.
$$
Then
\begin{equation}\label{eqn:MOtO}
\int_{[\overline{\Mor}(\Sigma_\ell,\Mbar(S,G^\C)]^\vir }\CO_1 \cdots \CO_k
=\int_{[\Mbar(\Sigma_\ell\times S,G^\C)]^\vir }\tilde{\CO}_1 \cdots \tilde{\CO}_k.
\end{equation}
The left and right hand sides of \eqref{eqn:MOtO} can be viewed as
the mathematical definition of the left and right hand sides of
\eqref{eqn:OtO}, respectively.

Here we propose some natural observables such that \eqref{eqn:MOtO} is
expected to hold. Given a point $p\in \Sigma_\ell$, define
\begin{eqnarray*}
\ev_p:\overline{\Mor}(\Sigma_\ell,\Mbar(S,G^\C))&\longrightarrow&
\Mbar(S,G^\C)\\
u &\mapsto & u(p)\\
r_{p}:\Mbar(\Sigma_\ell\times S,G^\C) &\longrightarrow&
\Mbar(S,G^\C)\\
P &\mapsto & P|_{\{p\}\times S}
\end{eqnarray*}
Then $\ev_p=r_p\circ \Phi$.
Given $\alpha_1,\ldots,\alpha_k\in H^*(\Mbar(S,G^\C))$ and $k$ distinct points
$p_1,\ldots, p_k \in \Sigma_\ell$, define
\begin{eqnarray*}
\langle\alpha_1\cdots\alpha_k\rangle_{\sigma|ASD}^{\Sigma_\ell | S}
&=&\int_{[\overline{\Mor}(\Sigma_\ell,\Mbar(S,G^\C)]^\vir }\ev_{p_1}^*\alpha_1\cup\cdots\cup \ev_{p_k}^* \alpha_k\\
\langle \alpha_1 \cdots\alpha_k \rangle_{DT}^{\Sigma_\ell\times S}
&=&\int_{[\Mbar(\Sigma_\ell\times S,G^\C)]^\vir }r_{p_1}^*\alpha_1\cup\cdots\cup r_{p_k}^* \alpha_k
\end{eqnarray*}
The definition is independent of choices of $p_1,\ldots, p_k$.
The correlation functions $\langle \alpha_1 \cdots \alpha_k \rangle_{\sigma |ASD}^{\Sigma_\ell |S}$
are $k$-point {\em mixed}\footnote{since we fix the complex structure of the domain} Gromov-Witten
invariants of $\Mbar(S,G^\C)$. We should have
\begin{equation}\label{eqn:ASD-DT}
\langle \alpha_1 \cdots \alpha_k\rangle_{\sigma| ASD}^{\Sigma_\ell| S} =
\langle \alpha_1 \cdots \alpha_k\rangle_{DT}^{\Sigma_\ell\times S}
\end{equation}

\subsection{Abelian theory}
In the abelian theory,
$$
\CM_{ASD} (S) \cong \Hilb^k S
$$
which is a resolution of the symmetric product orbifold, $\Sym^k S$.
In this case, DT theory on $\Sigma_\ell \times S$ is related to
GW theory on $\Sigma_\ell \times S$ by the GW/DT correspondence
conjectured in \cite{MNOP1, MNOP2}. These lead to the equivalence
of the following theories:
\begin{itemize}
\item[S1]  Gromov-Witten theory on $\Sigma_\ell\times S$.
\item[S2]  Donaldson-Thomas theory on $\Sigma_\ell \times S$.
\item[S3]  Mixed Gromov-Witten invariants with domain $\Sigma_\ell$ and target $\Hilb^k S$.
\end{itemize}

We now give a more precise conjecture relating S2 and S3.
We first introduce correlation functions in S3.
Classical cohomology ring of $\Hilb^k S$ is determined first for $\C^2$ in \cite{LS1, V}
and then for K3 \cite{LS2}. The general case is treated in \cite{CG}.
As vector spaces, there is an isomorphism:
$$
\CF= \Sym^*\Bigl(t^{-1}\Q[t^{-1}]\otimes H^*(S;\Q)\Bigr)
\cong \bigoplus_{k\geq 0} H^*(\Hilb^k S;\Q).
$$
A basis of $\CF$ is given by {\em cohomology weighted partitions}
$$
\alpha_{-\mu_1}(\gamma_{i_1})\cdots \alpha_{-\mu_n}(\gamma_{i_n})\mid 0\rangle.
$$
where $(\mu_1,\mu_2,\ldots,\mu_\ell)$ is a partition, and
$\{\gamma_i\}$ is a basis of the cohomology {\em group} $H^*(S;\Q)$.

Assume that $3\ell-3+n \geq 0$, so that the moduli space
$\Mbar_{\ell,n}$ of genus $\ell$, $n$-pointed stable curves is nonempty.
We fix $n$ distinct points $x_1,\ldots,x_n$ on $\Sigma_\ell$, and let
$\xi\in H^{6\ell-6+2n}(\Mbar_{\ell,n};\Q)$ be the
Poincar\'e dual of the point class $[(\Sigma_\ell, x_1,\ldots,x_n)] \in \Mbar_{\ell,n}$.
Given $\beta\in H_2(\Hilb^k S,\Z)$, let
$\pi: \Mbar_{\ell,n}(\Hilb^k S, \beta) \to \Mbar_{\ell,n}$ be the forgetful  map.
We have the following genus $\ell$, $n$-point, mixed Gromov-Witten invariants:
$$
\langle (\mu^1,\delta^1),\ldots,(\mu^n,\delta^n)\rangle_{\xi,\beta}
=\int_{[\Mbar_{\ell,n}(\Hilb^k S, \beta)]^{\vir}}\pi^*\xi\cup \prod_{i=1}^n \ev_i^* (\mu^i,\delta^i)
$$
where
$$
\mu^i\in \CP,  \quad |\mu^i|=k, \quad \delta^i \in H^*(S;\Q)^{\otimes \ell(\mu^i)}.
$$

Let $\{D_1,\ldots, D_m\}$ be a basis of $H^2(S;\Z)$. Then
$$
\tD_0= \alpha_{-2}(1)\left(\alpha_{-1}(1)\right)^{k-1} \mid 0\rangle ,\quad
\tD_i=\alpha_{-1}(D_i)\left(\alpha_{-1}(1)\right)^{k-1} \mid 0\rangle,\ i=1,\ldots,m
$$
form a basis of $H^2(\Hilb^k(S);\Z)$. Define
\begin{align}\label{eqn:Hilb}
&\langle (\mu^1,\delta^1),\ldots,(\mu^k,\delta^n)\rangle^{\Hilb^k(S)}_\xi
=\notag\\
&
\sum_{\beta\in H_2(\Hilb^k(S);\Z)} \sum_\beta \langle (\mu^1,\delta^1),\ldots,(\mu^k,\delta^n)\rangle_{\xi,\beta}
\tq_0^{\beta\cdot \tD_0}\prod_{j=1}^m \tq_j^{\beta\cdot \tD_j}
\end{align}

We next introduce correlation functions in S2.
We fix $n$ distince points $x_1,\ldots,x_n$ on $\Sigma_\ell$ as before. Given $\beta\in H_2(S;\Z)$, there
are maps
$$
\epsilon_i: I_s(\Sigma_\ell\times S/ \{x_1,\ldots,x_n\}\times S, k[\Sigma_\ell]+\beta) \to \Hilb^k(\{x_i\}\times S)\cong
\Hilb^k(S,\Z),\quad
i=1,\ldots, n.
$$
Let $(\mu^i,\delta^i)$ be cohomology weighted partitions as above. Define
$$
\langle(\mu^1,\delta^1),\ldots,(\mu^n,\delta^n)\rangle^{DT(\xi,S)}_{s,\beta}
=\int_{[I_s(\Sigma_\ell\times S/ \{x_1,\ldots,x_n\}\times S, k[\Sigma_\ell]+\beta)]^\vir}
\prod_{i=1}^n \epsilon_i^*(\mu^i,\delta^i).
$$
\begin{align}\label{eqn:relativeDT}
&\langle(\mu^1,\delta^1),\ldots,(\mu^n,\delta^n)\rangle^{DT(\xi,S)}_k
\notag\\
&
=q^{-k}\sum_{s\in \Z} \sum_{\beta\in H_2(S,\Z)}
\langle(\mu^1,\delta^1),\ldots,(\mu^n,\delta^n)\rangle^{DT(\xi,S)}_{s,\beta}
q^s \prod_{j=1}^m q_j ^{\beta\cdot D_j}
\end{align}

\begin{conj}[Maulik-Oblomkov]\label{conj:hyperkahler}
Let $S$ be a  hyper-K\"{a}hler surface. Then
\begin{equation}\label{eqn:MO}
\langle (\mu^1,\delta^1),\ldots,(\mu^n,\delta^n)\rangle^{\Hilb^k(S)}_\xi
=\langle(\mu^1,\delta^1),\ldots,(\mu^n,\delta^n)\rangle_k^{DT(\xi,S)}
\end{equation}
\end{conj}
Maulik and Oblomkov's original conjecture concerns the genus zero case
(i.e. $\Sigma_\ell= \cP^1$), but the general case should follow
from the genus 0 case by degenerating $\Sigma_\ell$ to
a union of $\cP^1$'s. Maulik and Oblomkov have checked that
Conjecture \ref{conj:hyperkahler} fails for Enrique surfaces.

\begin{conj} [Maulik-Oblomkov]\label{conj:Fano}
Let $S$ be a Fano surface, and let $\Sigma_\ell=\cP^1$.
Let $c_1(S) =\sum_{j=1}^m c_j D_j$. Then
\begin{equation}\label{eqn:MO-Fano}
\langle (\mu^1,\delta^1),\ldots,(\mu^n,\delta^n)\rangle^{\Hilb^k(S)}_\xi
=\langle(\mu^1,\delta^1),\ldots,(\mu^n,\delta^n)\rangle_k^{DT(\xi,S)}
\end{equation}
with $\tq_0 = q$, $\tq_j = (1-q^{-1})^{c_j} q_j$.
\end{conj}

Maulik and Oblomkov have checked Conjecture \ref{conj:Fano} for
$\cP^1\times \C$ and $\CO_{\cP^1}(-1)\to \cP^1$.

We summarize some known results which motivated
the above conjectures.
When the genus $\ell=0$, then 3-point correlation
functions in $Z_{2D}(\Sigma_0)=Z_{2D}(\cP^1)$ are the structure
constants of the quantum cohomlogy ring of $\CM_{ASD}(S)$.
Here we consider small quantum cohomolgy, so the structure constants
of the quantum product involves $3$-point genus
zero Gromov-Witten invariants. The big quantum cohomlogy involves
$n$-point genus zero Gromov-Witten invariants, which are not
the same as the $n$-point mixed Gromov-Witten invariants.
The crepant resolution conjecture
asserts that the quantum cohomology ring of $\Hilb^k S$ is isomorphic
to the orbifold quantum cohomology ring of $\Sym^k S$ \cite{Za, Va,CR}

When $S=\C^2$, $\Sigma_0\times S=\cP^1\times \C^2$ is non-compact, so
a priori neither GW theory nor DT theory is defined. Using
the torus action on $\C^2$, one can define equivariant GW
and DT theories.
\begin{enumerate}
\item[C1] Equivariant Gromov-Witten theory on $\cP^1\times \C^2$.
\item[C2] Equivariant Donaldson-Thomas theory on $\cP^1\times \C^2$.
\item[C3] Equivariant quantum cohomology of $\Hilb_k\C^2$.
\item[C4] Equivariant orbifold quantum cohomology of $\Sym^k\C^2$.
\end{enumerate}
The equivalences among C1, C2, C3 are proved by Bryan, Okounkov,
Pandharipande \cite{BP, OP1,OP2}, whereas
the equivalence between C3 and C4 is studied in \cite{BG}.

Maulik and Oblomkov established the equivalance of the following
theories when $S$ is an ADE resolution.
\begin{enumerate}
\item[A1] Equivariant Gromov-Witten theory on $\cP^1\times S$.
\item[A2] Equivariant Donaldson-Thomas theory on $\cP^1\times S$.
\item[A3] Equivariant quantum cohomology of $\Hilb_k S$.
\end{enumerate}

In particular, C2$\Leftrightarrow$C3 and A2$\Leftrightarrow$A3
are special cases of the equivariant version of \eqref{eqn:MO}.

We may consider the case where $S$ is a orbifold. 
\begin{enumerate}
\item[O1] Equivariant orbifold Gromov-Witten theory on $\cP^1\times [\C^2/\Z_{n+1}]$.
\item[O2] Equivariant Donaldson Thomas theory on $\cP^1\times[\C^2/\Z_{n+1}]$.
\item[O3] Equivariant quantum cohomology of $\Hilb_k[\C^2/\Z_{n+1}]$.
\end{enumerate}
By Maulik-Okounkov \cite{MO}, when $S=\mathcal{A}_n$ (resolution of $A_n$-singularity), A3 is
equivalent of O3, which can be viewed
as an example of Ruan's Crepant Tranformation Conjecture. Z. Zhou  \cite{Zhou} established
the equivalence between O2 and O3.

Finally, let $S$ be a smooth projective K3 surface. The Gromov-Witten theories of 
$K3$ and $K3\times \cP^1$ are trivial, so one considers reduced Gromov-Witten theory
defined by reduced virtula class. The compact case is much more difficult than the 
non-compact case. 
\begin{enumerate}
\item[R1] Reduced Gromov-Witten theory of $\cP^1\times S$
\item[R2] Reduced Donaldson-Thomas theory on $\cP^1\times S$.
\item[R3] Reduced quantum cohomology of $\Hilb_kS$.
\end{enumerate}
G. Oberdieck studies R1 and R3 in \cite{Ob2} and \cite{Ob1}, respectively. 
In \cite{Ob1}, Oberkieck conjectures a fomrula for the quantum multiplication
with divisor classes on $\Hilb_kS$, and prove the conjecture in the first
non-trivial case $\Hilb_2 S$.

\section{Dimensional reduction of torsion DT theory}
\subsection{Seiberg-Witten invariants via Mochizuki's wall crossing}\label{subsec:review}
Let $S$ be a smooth projective surface and $h$ an ample
divisor on it.
We assume that $H^1(\oO_S)=0$ and the arithmetic genus $p_g=\dim H^2(\oO_S)>0$.
(For instance, any smooth hypersurface in a quintic 3-fold
satisfies this assumption.)
Let us take an element
\begin{align*}
v=(2, a, n) \in H^0(S) \oplus H^2(S) \oplus H^4(S)
\end{align*}
such that $h \cdot a$ is an odd number.
Then by this choice, any $h$-semistable sheaf $E \in \Coh(S)$
with $\ch(E)=v$ is $h$-stable.
Now define $\mM_{h}(v)$ to be the moduli
space of $h$-stable sheaves $E \in \Coh(S)$
with $\ch(E)=v$.
For simplicity, we assume that
there exists a universal sheaf
$\eE \in \Coh(S \times \mM_h(v))$.
Let
$p_M$ be the projection from $S \times \mM_h(v)$
to $M_h(v)$. We have the
decomposition
\begin{align*}
\dR p_{M\ast} \dR \hH om(\eE, \eE) =\dR p_{M\ast} \dR \hH om(\eE, \eE)_0
\oplus \dR p_{M\ast} \oO_{S \times \mM_h(v)}
\end{align*}
and the perfect obstruction theory
\begin{align*}
E^{\bullet}_{\mM_{h}(v)}:=\dR p_{M\ast}
\dR \hH om(\eE, \eE)_0^{\vee} \to \dL_{\mM_{h}(v)}.
\end{align*}
We have the associated virtual cycle $[\mM_h(v)]^{\rm{vir}}$
whose virtual dimension $d$ is
\begin{align*}
d=a^2 -4n -3\chi(\oO_S).
\end{align*}
Now let $P(\eE)$ be a polynomial
in the slant products $\ch_i(\eE)/b$ for the
elements $b \in H^{\ast}(S)$ and $i\in \mathbb{Z}_{\ge 0}$.
By the wall crossing argument using the master space,
Mochizuki described the invariant
\begin{align*}
\int_{[M_h(v)]^{\rm{vir}}} P(\eE)
\end{align*}
in terms of the Seiberg-Witten invariants and
certain integration over the Hilbert schemes of points on $S$.

The SW invariant is then defined as follows:
for any $c \in \mathrm{NS}(S)$, let $L$ be the line
bundle on $S$ such that $c_1(L)=c$, which is uniquely determined following
the assumption that $H^1(\oO_S)=0$.
Let $M(c)$ be the moduli space of non-zero
morphisms $\oO_S \to L$, which is isomorphic to
$\mathbb{P}(H^0(S, L))$. The natural deformation theory
of pairs $\oO_S \to L$ induces an obstruction bundle
$O(c)$ on $M(c)$, which fits into the exact sequence
\begin{align*}
0 \to H^1(S, L) \otimes \oO_{M(c)}(1) \to O(c) &\to H^2(S, \oO_S) \otimes
\oO_{M(c)} \\
 &\to H^2(S, L) \otimes \oO_{M(c)}(1) \to 0.
\end{align*}
The virtual cycle $[M(c)]^{\rm{vir}}$ is defined to be
the Euler class of $O(c)$. If it is non-zero, then the
virtual dimension is zero\footnote{The $p_g>0$ assumption
is required here in~\cite[Proposition~6.3.1]{Moc}.}, and
\begin{align*}
\mathrm{SW}(c) = \int_{[M(c)]^{\rm{vir}}} 1.
\end{align*}
By setting $d(c)=h^0(S, L)-1$, the SW invariant
is computed as (cf.~\cite[Proposition~6.3.1]{Moc})
\begin{align*}
\mathrm{SW}(c)=(-1)^{d(c)}\left( \begin{array}{c}
p_g -1 \\
d(c)
\end{array}  \right).
\end{align*}

Let $S^{[n]}$ be the Hilbert scheme of $n$-points in $S$.
For $j=1, 2$, let
$\zZ_i \subset S \times S^{[n_i]}$ be the universal
subscheme and
$\iI_i \subset \oO_{S \times S^{[n_i]}}$ be
its ideal sheaf.
Below we consider the decomposition
\begin{align*}
a_1+a_2=a \,, \qquad a_i \in \NS(S).
\end{align*}
We denote by $e^{a_i}$ the
line bundle on $S$
whose first Chern class equals to $a_i$.
We define $Q(\iI_1 e^{a_1 -s}, \iI_2 e^{a_2 +s})$ to be the
Euler class of the following $\mathbb{C}^{\ast}$-equivariant
virtual vector bundle on $S^{[n_1]} \times S^{[n_2]}$:
\begin{align*}
-\dR p_{\ast} \dR \hH om(\iI_1 e^{a_1 -s}, \iI_2 e^{a_2+s})
-\dR p_{\ast} \dR \hH om(\iI_2 e^{a_2+s}, \iI_1 e^{a_1 -s}).
\end{align*}
Here $s$ is the equivariant parameter with respect to the
trivial $\mathbb{C}^{\ast}$-action,
and $p$ is the projection from $S \times S^{[n_1]} \times S^{[n_2]}$
to $S^{[n_1]} \times S^{[n_2]}$.
All the equivariant sheaves in the derived inner Hom's
are pulled back to $S \times S^{[n_1]} \times S^{[n_2]}$.

We also consider the rank $n_i$ vector bundle
on $S^{[n_i]}$, given by
\begin{align*}
\vV_i =
p_{\ast} (\oO_{\zZ_i} \otimes e^{a_i})
\end{align*}
We define $\aA(a_1, a_2, v)$ to be
\begin{align}\notag
&\aA(a_1, a_2, v) = \\
\label{def:A}
& \sum_{\begin{subarray}{c}
n_1+n_2=n-a_1 a_2 \\
n_1>n_2
\end{subarray}}
\int_{S^{[n_1]} \times S^{[n_2]}} \mathrm{Res}_{s=0}
\left( \frac{P(\iI_1 e^{a_1 -s} \oplus \iI_2 e^{a_2 +s})}{Q(\iI_1 e^{a_1 -s},
\iI_2 e^{a_2 +s})}  \frac{e(\vV_1) \cdot e(\vV_2 e^{2s})}{(2s)^{n_1 + n_2 -p_g}}\right).
\end{align}
The following result was obtained by Mochizuki:
\begin{theorem}\emph{(Mochizuki~\cite[Theorem~1.4.6]{Moc})}\label{thm:Moc}
Assume that $ah >2K_S h$ and $\chi(v) = \int_S v \cdot \td_S \ge 1$.
Then we have the following formula:
\begin{align*}
\frac{1}{2}\int_{\mM_h(v)}P(\eE) =-\sum_{\begin{subarray}{c}
a_1 + a_2 =a \\
a_1 h < a_2 h
\end{subarray}}
\mathrm{SW}(a_1) \cdot 2^{1-\chi(v)} \cdot \aA(a_1, a_2, v).
\end{align*}
\end{theorem}
\begin{remark}
Note that the factor $1/2$ in the LHS
comes from the difference between Mochizuki's convention and ours.
Mochizuki used the moduli space of oriented stable sheaves,
which is a $\mu_2$-gerb over our moduli space $M_h(v)$.
\end{remark}
\begin{remark}
The assumptions $ah >2K_S h$ and $\chi(v) \ge 1$
are satisfied if we replace $v$ by $v \cdot e^{kh}$ for $k\gg 0$.
\end{remark}

\subsection{Torsion DT invariants and Seiberg-Witten  theory}
Let $(S, h)$ and $v\in H^{\ast}(S, \mathbb{Q})$
be as in the previous subsection, and consider
\begin{align}
\pi \colon X_{S} = \omega_S \to S
\end{align}
the total space of the canonical line bundle on $S$.
Note that $X_{S}$ is a non-compact Calabi-Yau 3-fold, i.e.
$\omega_{X_{S}} \cong \oO_{X_{S}}$.
Denote
\begin{align}
\Coh_c(X_{S}) \subset \Coh(X_{S})
\end{align}
to be the abelian category of
coherent sheaves on $X_{S}$
whose supports are compact.
The slope function $\mu_h$ on $\Coh_c(X_{S}) \setminus \{0\}$ defined as
\begin{align*}
\mu_h(E) =\frac{c_1(\pi_{\ast}E) \cdot h}{\rank(\pi_{\ast}E)} \in
\mathbb{Q} \cup \{ \infty\}
\end{align*}
determines a slope stability condition on $\Coh_c(X_{S})$ in
the usual way.
Let $\overline{\mM}_h(v)$ be the moduli space
of $\mu_h$-stable sheaves $E \in \Coh_c(X_{S})$ with $\ch(\pi_{\ast}E)=v=<r, \gamma, n>$.
The DT invariant for the local surface $X$ is then defined by
\begin{align}\label{DThv}
\DT_h(v)=\int_{\overline{\mM}_h(v)} \nu_M \,\, d\chi.
\end{align}
Here $\nu_M$ is Behrend's constructible function~\cite{Beh} on
$\overline{\mM}_h(v)$.

The one-dimensional complex torus $\mathbb{C}^{\ast}$
acts on $X_{S}$ by re-scaling on the fibers of $\pi$.
By the localization, the DT invariant (\ref{DThv})
coincides with the integration of $\nu_M$
over the $\mathbb{C}^{\ast}$-fixed locus
$\overline{\mM}_h(v)^{\mathbb{C}^{\ast}}$.
Note that the moduli space $\mM_h(v)$
is an open and closed subscheme of
$\overline{\mM}_h(v)^{\mathbb{C}^{\ast}}$.

Now let $\overline{\eE} \in \Coh(X_{S} \times \overline{\mM}_h(v))$
be the universal family,
and $p_{\overline{M}}$ the projection
from $X_{S} \times \overline{\mM}_h(v)$
to $\overline{\mM}_h(v)$.
Let $\dR \hH om(\overline{\eE}, \overline{\eE})_0$
be the cone of the composition
\begin{align*}
\dR \hH om_{\overline{\mM}_h(v) \times X_{S}}
(\overline{\eE}, \overline{\eE}) \stackrel{\pi_{\ast}}{\to}
\dR \hH om_{\overline{\mM}_h(v) \times S}(\pi_{\ast} \overline{\eE},
\pi_{\ast}\overline{\eE}) \stackrel{\mathrm{tr}}{\to}
\oO_{\overline{\mM}_h(v) \times S}.
\end{align*}
We obtain the $\mathbb{C}^{\ast}$-fixed trace-free
perfect\footnote{By the Serre duality and the non-trivial
$\mathbb{C}^{\ast}$-weight on
$\omega_S$, the higher obstruction space vanishes
after taking the $\mathbb{C}^{\ast}$-fixed part.}
obstruction theory
\begin{align}\label{3fold}
E^{\bullet}_{\overline{\mM}_h(v)}:=\left( \dR p_{\overline{M}\ast} \dR \hH om(\overline{\eE},
\overline{\eE})_0[1]
\right)^{\vee\mathbb{C}^{\ast}}
\to \dL_{\overline{\mM}_h(v)^{\mathbb{C}^{\ast}}}.
\end{align}
Let $[\overline{\mM}_h(v)^{\mathbb{C}^{\ast}}]^{\rm{vir}}$
be the associated virtual fundamental class.
Instead of working with the invariant (\ref{DThv}),
we can consider the invariant
\begin{align} \label{DThv3}
\overline{\overline{\DT}}_h(v) =
\int_{[\overline{\mM}_h(v)^{\mathbb{C}^{\ast}}]^{\rm{vir}}}
c\left( \left( \dR p_{\overline{M}\ast} \dR \hH om(\overline{\eE}, \overline{\eE})_0[1]
\right)^{\vee\mathbb{C}^{\ast}}\right).
\end{align}
\begin{remark}
Here the definition of the invariant (\ref{DThv3})
is motivated by Fantechi-G\"ottsche's virtual
Euler numbers~\cite[Section~4]{FG}. In fact, it coincides with the virtual
Euler number of $\overline{\mM}_h(v)$ up to a sign.
In particular,
when $\overline{\mM}_h(v)^{\mathbb{C}^{\ast}}$ is non-singular with
the expected dimension, then the above invariant coincides with the
topological Euler number of $\overline{\mM}_h(v)$ up to a sign.
\end{remark}

Now let $r=2$. The $\mathbb{C}^{\ast}$-fixed locus
$\overline{\mM}_h(v)^{\mathbb{C}^{\ast}}$
decomposes into two components
\begin{align*}
\overline{\mM}_h(v)^{\mathbb{C}^{\ast}}=
\mM_h(v) \coprod \widehat{\mM}_h(v)
\end{align*}
where $\mM_h(v)$ is, roughly speaking, the moduli space of $\mu_h$-stable torsion-free sheaves of rank $2$ on $S$, and
$\widehat{\mM}_h(v)$ is the moduli space of
$\mathbb{C}^{\ast}$-fixed
$\mu_h$-stable sheaves with rank one
on the ``fat'' surface $2S$. Now we can define the DT invariants associated to each component, by restricting the obstruction theory in \eqref{3fold} to each component and integrating against the corresponding  induced virtual cycle. Let
\begin{align}\label{DThv4}
\widehat{\DT}_h(v) =
\int_{[\widehat{\mM}_h(v)]^{\rm{vir}}}
c\left( \left( \dR p_{\overline{M}\ast} \dR \hH om(\overline{\eE}, \overline{\eE})_0[1]
\right)^{\vee\mathbb{C}^{\ast}} \right)
\end{align}
be the contribution from the ``fat'' surface $2S$ and
\begin{align}\label{DThv4}
\DT_h(v) =
\int_{[\mM_h(v)]^{\rm{vir}}}
c\left( \left( \dR p_{\overline{M}\ast} \dR \hH om(\overline{\eE}, \overline{\eE})_0[1]
\right)^{\vee\mathbb{C}^{\ast}} \right)
\end{align}
 be the contribution from rank $2$ sheaves on $S$. Note that, as we saw in the last section, the latter contributions have been computed in some cases by Mochizuki \cite{Moc} (look at Theorem \ref{thm:Moc}). On the other hand, it can be shown \cite[Proposition 3.13, 3.14, Corollary 3.15, and Proposition 3.21]{GLSY} that the contribution of $\widehat{\mM}_h(v)$ to $\widehat{\DT}_h(v)$ is given by the invariants of ``\textit{nested Hilbert schemes"} on $S$, denoted by $S^{[n_{1}>n_{2}]}_{\beta}$, parametrizing 2-step flags $\iI_{Z_{1}}(-C_{1})\hookrightarrow \iI_{Z_{2}}$ of (possibly twisted) ideal sheaves of subschemes, $(Z_{1}, Z_{2}), (C_{1})$ of $S$, where $Z_{i}, i=1,2$ are zero dimensional subschemes of $S$ with length $n_1,n_2$ respectively, $C_{1}\subset S$ is a divisor with $[C_{1}]=\beta$ for some suitable $\beta$, such that $Z_{2}$ is a subscheme of $Z_{1}\cup C_{1}$ (hence it induces the injective map of the corresponding ideal sheaves). We then prove the following identity in \cite[Theorem 7]{GLSY} relating $\overline{\DT}_h(v)$ to SW invariants and the invariants of nested Hilbert schemes contributing to $\widehat{\DT}_h(v)$.
\begin{equation}\label{SW/DT}
\overline{\overline{\DT}}_h(v)=-\sum_{\begin{subarray}{c}
a_1 + a_2 =a \\
a_1 h < a_2 h
\end{subarray}}
\mathrm{SW}(a_1) \cdot 2^{1-\chi(v)} \cdot \aA(a_1, a_2, v)+\widehat{\DT}_h(v).
\end{equation}
Here $\widehat{\DT}_h(v)$ is, roughly speaking, given as a sum over the contribution of nested Hilbert schemes $S^{[n_1>n_2]}_{\beta}$, for all allowed values of $n_{1}, n_{2}, \beta$ induced by the choice of $v$. It must be pointed out that by S-duality consideration, the generating series of the invariants on the left-hand side of \eqref{SW/DT} is expected to be given by a vector valued modular form of weight $-\frac{3}{2}$, while the generating series of $\widehat{\DT}_h(v)$ is also (in some cases) shown to be given by modular forms of certain weight \cite[Theorem 5]{GLSY}. Therefore, equation \eqref{SW/DT} implies that the generating series of the Seiberg-Witten invariants (despite not being necessarily modular itself) can in some cases be written with respect to modular forms.



\begin{thebibliography}{99}

\bibitem{AB} B.~Anchouche, I.~Biswas,
``Einstein-Hermitian connections on polystable principa bundles
over a compact K\"{a}hler manifolds'', Amer. J. Math. {\bf 123} (2001),
207--228.

\bibitem{AKMV} M.~Aganagic, A.~Klemm, M.~Mari\~{n}o, C.~Vafa,
``The topological vertex'', Commun. Math. Phys. {\bf 254} (2005)
425--478.

\bibitem{ALSpence} B.~S.~Acharya, M.~O'Loughlin and B.~J.~Spence,
``Higher-dimensional analogues of Donaldson-Witten theory,''
Nucl.\ Phys.\ B {\bf 503} (1997) 657, hep-th/9705138.


\bibitem{AMV} M. Aganagic, M. Mari\~{n}o, C. Vafa,
``All loop topological string Amplitudes from Chern-Simons theory,'' 
Comm, Math. Phys., {\bf 247}, (2004) no 2, 467--512

\bibitem{AOSV} M.~Aganagic, H.~Ooguri, N.~Saulina, C.~Vafa,
``Black Holes, q-deformed 2d Yang-Mills and nonperturbative topological strings'', 
Nucl. Phys. B.
{\bf 715}, (2005) no 1�2, 23, 304--348.

\bibitem{BKSinger} L.~Baulieu, H.~Kanno and I.~M.~Singer,
``Special quantum field theories in eight and other dimensions,''
Commun.\ Math.\ Phys.\  {\bf 194} (1998), no.1,  149--175.

\bibitem{Beh} K.~Behrend,
``Donaldson-Thomas invariants via microlocal geometry,''
Ann. of Math. (2) {\bf 170} (2009), no. 3, 1307--1338.


\bibitem{BCOV} M.~Bershadsky, S.~Cecotti, H.~Ooguri, C.~Vafa,
``Holomorphic Anomalies in Topological Field Theories,''
Nuclear Phys. {\bf B 405} (1993), no. 2-3, 279-304.


\bibitem{BJSV} M. Bershadsky, A. Johansen, V. Sadov, C. Vafa,
``Topological Reduction of 4D SYM to 2D $\sigma$--Models,''
hep-th/9501096.



\bibitem{BSV} M.~Bershadsky, C.~Vafa and V.~Sadov,
``D-Branes and Topological Field Theories,''
Nucl.\ Phys.\ B {\bf 463} (1996) 420, hep-th/9511222.


\bibitem{BT} M.~Blau and G.~Thompson,
``Euclidean SYM theories by time reduction and special holonomy  manifolds,''
Phys.\ Lett.\ B {\bf 415} (1997) 242, hep-th/9706225.

\bibitem{BT} M.~Blau and G.~Thompson,
``Euclidean SYM theories by time reduction and special holonomy manifolds,''
Phys. Lett. {\bf B 415} (1997), no. 3, 242--252.

\bibitem{BG} J. Bryan, T. Graber,
``The crepant resolution conjecture,''
Algebraic geometry—Seattle 2005. Part 1, 23--42,
Proc. Sympos. Pure Math., {\bf 80}, Part 1, Amer. Math. Soc., Providence, RI, 2009.

\bibitem{BL}
J.~Bryan,  N.C.~Leung, 
``The enumerative geometry of $K3$ surfaces and modular forms,'' 
J. Amer. Math. Soc. {\bf 13} (2000), no. 2, 371--410.

\bibitem{BP} J.~Bryan, R.~Pandharipande,
``The local Gromov-Witten theory of curves,''
With an appendix by Bryan, C. Faber, A. Okounkov and Pandharipande.
J. Amer. Math. Soc. {\bf 21} (2008), no. 1, 101--136.



\bibitem{CKYZ} T.-M.~Chiang, A.~Klemm, S.-T.~Yau, E.~Zaslow,
``Local Mirror Symmetry: Calculations and Interpretations,''
hep-th/9903053, Adv. Theor. Math. Phys. {\bf 3}  (1999),
no. 3, 495--565.

\bibitem{CG}
K. Costello, I. Grojnowski,
``Hilbert schemes, Hecke algebras and the Calogero-Sutherland system,''
math.AG/0310189.

\bibitem{CR} W. Chen, Y. Ruan,
``Orbifold Gromov-Witten theory,'' 
{\em Orbifolds in mathematics and physics (Madison, WI, 2001)}, 25--85,
Contemp. Math., {\bf 310}, Amer. Math. Soc., Providence, RI, 2002.


\bibitem{Dij} R.~Dijkgraaf,
``Mirror symmetry and elliptic curves,''
The moduli space of curves (Texel Island, 1994), 149--163,
Progr. Math., {\bf 129}, Birkh\"{a}user Boston, Boston, MA, 1995.

\bibitem{DM} R.~Dijkgraaf and G.~W.~Moore,
``Balanced topological field theories,''
Commun. Math. Phys.  {\bf 185} (1997), no. 2, 411–-440.

\bibitem{DMVVsym}
R.~Dijkgraaf, G.~W.~Moore, E.~Verlinde and H.~Verlinde, ``Elliptic
genera of symmetric products and second quantized strings,''
Commun. Math. Phys.  {\bf 185} (1997)  no. 1, 197--209.

\bibitem{DF} D.-E. Diaconescu, B. Florea,
``The ruled vertex and nontoric del Pezzo surfaces,'' 
J. Hep {\bf 12}, (2006) no. 028.

\bibitem{DFS} D.-E.~Diaconescu, B.~Florea, N.~Saulina,
``A vertex formualism for local ruled surfaces,'' 
Comm. Math. Phys. {\bf 265}, (2006) no. 1, 201--226.


\bibitem{Do} S.~Donaldson,
``Anti-self-dual Yang-Mills connections over complex algebraic surfaces and stable vector bundles,'' 
Proc. Lond. Math. Soc. {\bf 50} (1985).

\bibitem{DT} S.~Donaldson and R.~Thomas, ``Gauge theory in higher dimensions,''
in {\it The geometric universe: science, geometry, and the work of Roger Penrose},
S. Huggett et. al eds., Oxford Univ. Press, 1998.

\bibitem{EQ} D. Edidin, Z. Qin,
``The Gromov-Witten and Donaldson-Thomas correspondence for trivial elleptic fibrations,''
Internat. J. Math. {\bf 18} (2007), no. 7, 821--838.

\bibitem{FG} B.~Fantechi and L.~G\"{o}ttsche,
``Riemann-Roch theorems and elliptic genus for virtually smooth schemes,''
Geom. Topol. {\bf 14} (2010), no. 1, 83--115.


\bibitem{FP} C.~Faber, R.~Pandharipande, 
``Hodge integrals and Gromov-Witten theory,'' 
Invent. Math. {\bf 139} (2000), no. 1, 173--199.

\bibitem{FS} R.~ Fintushel, R.~ Stern, 
``Immersed spheres in 4-manifolds and the immersed Thom conjecture,'' 
Turkish Journal of Mathematics {\bf 19} (1995), 145--157.

\bibitem{Gottsche} L.~G\"{o}ttsche, `
``The Betti numbers of the Hilbert scheme of points on a smooth projective surface,''  
Math. Ann. {\bf 286} (1990) 193--207.

\bibitem{GS} T. L. G\'{o}mez, I. Sols,
``Projective moduli space of semistable principal sheaves for a reductive group'',
Comm. Alg. Alg. Geom., Catania,  (2001).

\bibitem{GLSY} A.~Gholampour, A.~Sheshmani, S.T.~Yau,
``Nested Hilbert schemes and local Donaldson-Thomas theory,'' 
in preparation.

\bibitem{GP}
L.~G\"{o}ttsche, R.~Pandharipande, 
``The quantum cohomology of blow-ups of $P\sp 2$ and enumerative geometry,'' 
J. Differential Geom. {\bf 48} (1998), no. 1, 61--90.

\bibitem{GHarris}
P.~Griffiths, J.~Harris, 
{\em Principles of Algebraic Geometry},
Pure and Applied Mathematics. Wiley-Interscience [John Wiley \& Sons], 
New York, 1978. xii+813 pp.
 

\bibitem{GViii}
R.~Gopakumar and C.~Vafa,
``M-theory and topological strings. I,II,''
 hep-th/9809187; hep-th/9812127.

\bibitem{HHofer}
F.~Hirzebruch, T.~H\"ofer, ``On the Euler number of an orbifold,''
Math. Ann. {\bf 286} (1990) 255.

\bibitem{HP} C. Hofman, J-S.~Park,
``Cohomological Yang-Mills Theories on Kahler 3-folds,''
Nuclear Phys. B {\bf 600}  (2001),  no. 1, 133--162.

\bibitem{HSTa} S.~Hosono, M.-H.~Saito, A.~Takahashi, 
``Holomorphic Anomaly Equation and BPS State Counting of Rational Elliptic Surface,''
Adv. Theor. Math. Phys. {\bf 3} (1999) 177--208.

\bibitem{HSTb} S.~Hosono, M.-H.~Saito, A.~Takahashi, 
``Relative Lefschetz Action and BPS State Counting,'' 
Internat. Math. Res. Notices, (2001), No. 15, 783--816.

\bibitem{Ho} S.~Hosono, 
``Countin BPS States via Holomorphic Anomaly Equations,''
Calabi-Yau varieties and mirror symmetry (Toronto, ON, 2001), 57--86, 
Fields Inst. Commun., 38, Amer. Math. Soc., Providence, RI, 2003.

\bibitem{Hu}
J.~Hu, "Gromov-Witten invariants of blow-ups along surfaces",
 Compositio Math. {\bf 125} (2001), no. 3, 345--352.


\bibitem{IP} E.-N.~Ionel, T.~Parker,
``The Gromov invariants of Ruan-Tian and Taubes'',
Math. Res. Lett. {\bf 4} (1997), no. 4, 521--532.

\bibitem{Iq} A.~Iqbal,
``All Genus Topological String Amplitudes and 5-brane Webs as
Feynman Diagrams'', hep-th/0207114.


\bibitem{INOV} A.~Iqbal, N.~Nekrasov, A.~Okounkov and C.~Vafa,
``Quantum foam and topological strings,'' 
J. High Energy Phys. 2008, no. 4, 011, 47 pp. 

\bibitem{JS} D.~Joyce, Y.~Song,  
``A theory of generalized Donaldson-Thomas invariants,'' 
Mem.  Amer. Math. Soc., Publication {\bf 217}, (2012), no. 1020.

\bibitem{KMW} Albrecht Klemm, Jan Manschot, Thomas Wotschke
``Quantum geometry of elliptic Calabi-Yau manifolds,''
Commun. Number Theory Phys. {\bf 6} (2012), no. 4, 849--917.

\bibitem{KM} M.~Kontsevich, Y.~Manin,
 ``Gromov-Witten classes, quantum cohomology, and enumerative geometry,''
Comm. Math. Phys. {\bf 164}  (1994),  no. 3, 525--562.

\bibitem{KS} M.~Kontsevich, Y.~Soibelman, 
``Stability structures, motivic Donaldson-Thomas invariants and cluster transformations,'' 
math.AG/0811.2435.

\bibitem{KKV}
S. Katz, A. Klemm, C. Vafa,
``M-Theory, Topological Strings, and Spinning Black Holes'', 
Adv. Theor. Math. Phys. {\bf 3} (1999), 1445--1537.

\bibitem{KW}
A. Kapustin, E. Witten,
``Electric-Magnetic Duality and the Geometric Langlands Program,''
Commun. Number Theory Phys. {\bf 1} (2007), no. 1, 1--236. 

\bibitem{LabastidaL}  J.~M.~F.~Labastida and C.~Lozano,
``Mathai-Quillen formulation of twisted N=4 supersymmetric gauge theories in four-dimensions,''
Nuclear Phys. B {\bf 502} (1997), no. 3, 741--790.

\bibitem{LMruled}
C.~Lozano and M.~Mari\~{n}o, 
``Donaldson invariants of product ruled surfaces and two-dimensional  gauge theories,'' 
Commun.\ Math.\ Phys.\  {\bf 220} (2001) no. 2, 231-261. 

\bibitem{LS1} M.~Lehn, C.~Sorger,
``Symmetric groups and the cup product on the cohomology of Hilbert schemes,''
Duke Math. J., {\bf 110}, (2001) no 2, 345--357.

\bibitem{LS2}: M.~Lehn, C.~Sorger,
``The cup product of the Hilbert scheme for K3 surfaces,'' 
Invent. math.{\bf 152}, (2003) no. 2, 305--329.

\bibitem{LR} A.M. Li, and Y.B. Ruan,
``Symplectic surgery and Gromov-Witten invariants of Calabi-Yau 3-folds,''
Invent. Math. {\bf 145} (2001), no. 1, 151--218.

\bibitem{LLLZ} J. Li, C.-C. M. Liu, K. Liu, and J. Zhou,
``A mathematical theory of the topological vertex,'' 
Geom. Topol. {\bf 13} (2009), no. 1, 527--621.

\bibitem{Li} J. Li,
``A Degeneration Formula of GW-Invariants,''
J. Diff. Geom. {\bf 60} (2002), no. 2, 177--354.

\bibitem{LL} T.-J.~ Li, A.-K.~Liu, 
``The equivalence between SW and Gr in the case where $b^+=1$,'' 
Internat. Math. Res. Notices, {\bf 7} (1999) 335--345;
``Uniqueness of symplectic canonical class, surface cone and
symplectic cone of 4-manifolds with $b^+ = 1$,'' 
J. Diff. Geom. {\bf 58} (2001) 331--370.

\bibitem{LLZ1} C.-C. M. Liu, K. Liu, and J. Zhou,
``A proof of a conjecture of Mari\~{n}o-Vafa on Hodge integrals,''
J. Differential Geom. {\bf 65} (2003), no. 2, 289--340. 

\bibitem{LLZ2} C.-C. M. Liu, K. Liu, and J. Zhou, 
``A formula of two-partition Hodge integrals,''
J. Amer. Math. Soc. {\bf 20} (2007), no. 1, 149--184.

\bibitem{LMruled} C.~Lozano and M.~Mari\~{n}o,
``Donaldson invariants of product ruled surfaces and two-dimensional  gauge theories,''
Commun. Math. Phys. {\bf 220} (2001), no. 2, 231--261.

\bibitem{Marcus}  N.~Marcus,
``The Other topological twisting of N=4 Yang-Mills,''
Nucl.\ Phys.\ B {\bf 452} (1995), no. 1-2, 331--345. 

\bibitem{MMsimpl} M.~Mari\~{n}o and G.~W.~Moore, 
``Donaldson invariants for non-simply connected manifolds,'' 
Commun.\ Math.\ Phys.\  {\bf 203} (1999), no. 2, 249--267.
(1999) hep-th/9804104.

\bibitem{MNOP1} D. Maulik, N. Nekrasov, A. Okounkov, R. Pandharipande,
``Gromov-Witten theory and Donaldson-Thomas Theory I,'' 
Compos. Math. {\bf 142}, (2006) no. 5, 1263--1285.
.

\bibitem{MNOP2} D. Maulik, N. Nekrasov, A. Okounkov, R. Pandharipande,
``Gromov-Witten theory and Donaldson-Thomas Theory II,''
Compos. Math. {\bf 142}, (2006) no. 5, 1286--1304.


\bibitem{MOOP} D. Maulik, A. Oblomkov, A. Okounkov, R. Pandharipande, 
``Gromov-Witten/Donaldson-Thomas correspondence for toric 3-folds,''
Invent. math. {\bf 186}, (2011) no. 2, 435--479.

\bibitem{MO} D. Maulik, A. Okounkov,
``Quantum Groups and Quantum Cohomology," {\tt arXiv:1211.1287}.

\bibitem{MP} D. Maulik, R. Pandharipande,
``Gromov-Witten theory and Noether-Lefschetz theory,'' 
A celebration of algebraic geometry, 469--507, 
Clay Math. Proc., {\bf 18}, Amer. Math. Soc., Providence, RI, 2013.


\bibitem{MNVW} J.~A.~Minahan, D.~Nemeschansky, C.~Vafa and N.~P.~Warner,
``E-strings and N = 4 topological Yang-Mills theories,''
Nucl.\ Phys.\ B {\bf 527} (1998) no.3, 581-623.

\bibitem{Sa} K. Sakai,
``Topological string amplitudes for the local half K3 surface,''
arXiv:1111.3967.

\bibitem{Moc} T.~Mochizuki,
{\em Donaldson Type Invariants for Algebraic Surfaces},
Lect. Notes Math., Springer-Verlag, Berlin, Vol. 1972, (2009).

\bibitem{WMoore}
G.~Moore, E.~Witten, 
``Integration over the u-plane in Donaldson theory,''
Adv. Theor. Math. Phys. {\bf 1} (1998) no.2, 298--387.

\bibitem{Ob1} G. Oberdieck, 
``Gromov-Witten invariants of the Hilbert schemes of points of a K3 surface,'' 
arXiv:1406.1139.

\bibitem{Ob2} G. Oberdieck, 
``Gromov-Witten theory of K3 x $\mathbb{P}^1$ and quasi-Jacobi forms,''
arXiv:1605.05238.

\bibitem{ObP} G. Oberdieck, R. Pandharipande,
``Curve counting on $K3 \times E$, the Igusa cusp form $\chi_{10}$, and descendent integration,''
arXiv:1411.1514.

\bibitem{ObS} G. Oberdieck, J. Shen,
``Curve counting on elliptically fibered Calabi-Yau 3-folds,''
arXiv:1608.07073.

\bibitem{OkonekT} C.~Okonek and A.~Teleman,
``Seiberg-Witten invariants for manifolds with $b_+=1$, and the universal wall crossing formula,'' 
Internat. J. Math. {\bf 7} (1996), no. 6, 811--832. 

\bibitem{OP} A. Okounkov, R. Pandharipande, 
``Hodge integrals and invariants of the unknot,'' 
Geom. Topol. {\bf 8} (2004), 675--699. 

\bibitem{OP1} A. Okounkov, R. Pandharipande, 
``Quantum cohomology for the Hilbert scheme of points in the plane,''
Invent. math., {\bf 179}, (2010) no. 3, 523--557.

\bibitem{OP2} A. Okounkov, R. Pandharipande, 
``The local Donaldson-Thomas theory of curves,'' 
Geom. Topol., {\bf 14} (2010), 1503--1567.

\bibitem{R} Y. Ruan, 
``Symplectic Topology and Complex Surfaces”,
{\em Geometry and analysis on complex manifolds}, 171--197, World Sci. Publ., River Edge, NJ, 1994.

\bibitem{Sch} A. Schmitt, 
``Singular principal bundles over higher-dimsional manifolds and their moduli spaces,''
Int. Math. Res. Not. 2002, no. 23, 1183--1209.

\bibitem{Taubesbook} C.~Taubes, 
{\em Seiberg Witten and Gromov invariants for symplectic 4-manifolds},
International Press, Somerville, MA, 2000.

\bibitem{TY} G.~Tian, S.-T.~Yau,
``K\"{a}hler-Einstein metrics on complex surfaces with $c_1>0$,''
Comm. Math. Phys.  {\bf 112}  (1987),  no. 1, 175--203.

\bibitem{T} G.~Tian,
``On Calabi's conjecture for complex surfaces with positive first Chern class,'' 
Invent. Math. {\bf 101} (1990), no. 1, 101--172.

\bibitem{TY1}G.~Tian and S.-T.~Yau,
``Complete K\"{a}hler manifolds with zero Ricci curvature. I,''
J. Amer. Math. Soc. {\bf 3}  (1990),  no. 3, 579--609.

\bibitem{TY2} G.~Tian and S.-T.~Yau,
``Complete K\"{a}hler manifolds with zero Ricci curvature. II,''
Invent. Math. {\bf 106}  (1991),  no. 1, 27--60.


\bibitem{Th} R.P. Thomas, 
``A holomorphic Casson invariants for Calabi-Yau 3-folds, and bundles on $K3$ fibrations,'' 
J. Diff. Geom. {\bf 54}, (2000) no. 2, 367--438

\bibitem{UY} K.~Uhlenbeck and S.-T. Yau, 
``On the existence of Hermitian-Yang-Mills connections in stable vector bundles'', 
Comm. Pure Appl. Math. {\bf 39} (1986), no. S, supp. S257--S293. 

\bibitem{Va} C. Vafa,
``String vacua and orbifoldized LG models'',
Modern Phys. Lett. A {\bf 4} (1989), no. 12, 1169--1185.

\bibitem{VW} C.~Vafa and E.~Witten,
``A Strong coupling test of S duality,''
Nucl.\ Phys.\ B {\bf 431} (1994), no. 1-2, 3--77.   

\bibitem{V} E.~Vasserot,
``Sur l'anneau de cohomologie du sch\'{e}ma de Hilbert de $\C^2$,''
C. R. Acad. Sci. Paris, S\'{e}. I Math. {\bf 332} (2001), no. 1, 7--12. 

\bibitem{Wmonopoles} E.~Witten, 
``Monopoles and four-manifolds,'' 
Math.\ Res.\ Lett.\ {\bf 1} (1994), no. 6, 769--796. 

\bibitem{WKahler} E.~Witten, 
``Supersymmetric Yang-Mills theory on a four manifold,''
J.\ Math.\ Phys.\  {\bf 35} (1994), no. 10, 5101--5135. 

\bibitem{Y} S.-T.~Yau,
``On the Ricci curvature of a compact K\"{a}hler manifold and the complex Monge-Amp\`{e}re equation. I,''
Comm. Pure Appl. Math. {\bf 31} (1978), no. 3, 339--411.

\bibitem{YZ} S.-T.~Yau, E.~Zaslow,
``BPS states, string duality, and nodal curves on K3,''
Nuclear Phys. {\bf B 471} (1996), no. 3, 503--512.

\bibitem{Za} E. Zaslow, 
``Topological orbifold models and quantum cohomology rings,'' 
Comm. Math. Phys. {\bf 156} (1993), no. 2, 301--331.

\bibitem{Zhou} Z. Zhou,
``Donaldson-Thomas theory of $[\mathbb{C}^2/\mathbb{Z}_{n+1}]\times \mathbb{P}^1$,"
arXiv:1510.00871.

\end{thebibliography}
\end{document}